\documentclass[11pt]{article}

\usepackage{fullpage}
\usepackage[letterpaper, margin=1in]{geometry}
\usepackage{setspace}
\usepackage{amsmath}
\usepackage{siunitx}
\usepackage{mathrsfs} 
\usepackage{epsfig}
\usepackage{graphicx}
\usepackage{subcaption}
\usepackage{latexsym}
\usepackage{epstopdf}
\usepackage{amssymb}
\usepackage{booktabs}
\usepackage{threeparttable}
\usepackage{color}
\usepackage{apacite}
\usepackage{multirow}
\usepackage{adjustbox}
\usepackage[hidelinks]{hyperref}
\usepackage{natbib}
\usepackage{float}
\usepackage{siunitx}
\usepackage{booktabs}
\usepackage{tocbibind}
\usepackage[font=small,labelfont=bf,labelsep=none]{caption}
\usepackage{makecell}
\usepackage{tabularx}
\makeatletter  
\newif\if@restonecol  
\makeatother

\usepackage[linesnumbered,ruled,vlined]{algorithm2e}
\usepackage{algpseudocode}

\usepackage{tabularx}
\usepackage{changes}
\usepackage{tikz}
\usetikzlibrary {arrows.meta}

\usepackage{enumitem}
\setenumerate[1]{itemsep=0pt,partopsep=0pt,parsep=0pt,topsep=0pt}
\setitemize[1]{itemsep=0pt,partopsep=0pt,parsep=0pt,topsep=0pt}
\setdescription{itemsep=0pt,partopsep=0pt,parsep=0pt,topsep=0pt}

\begin{document}

\newtheorem{theorem}{Theorem}
\newtheorem{prop}{Proposition}
\newtheorem{lemma}{Lemma}
\newtheorem{remark}{Remark}
\newtheorem{corollary}{Corollary}
\newtheorem{observation}{Observation}
\newtheorem{property}{Property}
\newtheorem{claim}{Claim}
\newtheorem{conj}{Conjecture}
\newtheorem{definition}{Definition}
\newtheorem{class}{Class}
\def\square{\vrule height6pt width7pt depth1pt}
\def\endpf{\hfill\square\bigskip}
\makeatletter
\newcommand{\rmnum}[1]{\romannumeral #1}
\newcommand{\Rmnum}[1]{\expandafter\@slowromancap\romannumeral #1@}
\makeatother
\bibliographystyle{apacite}
\sisetup{group-digits = true}

\begin{center}
{\bf\large A $1.431$-Competitive Algorithm for Combinatorial Group Testing}\\\vspace{0.5cm}
Jun Wu$^{\, a, b}$, Yongxi Cheng$^{\, a, c}$, 
Zhen Yang$^{\, a, c,*}$
, Feng Chu$^{\,b}$, Junkai He$^{\, d}$
\\\vspace{3mm}
{\small
$^a$School of Management, Xi'an Jiaotong University, Xi'an 710049, China\\
{\footnotesize E-mail addresses: wjun048@stu.xjtu.edu.cn, chengyx@mail.xjtu.edu.cn, zhen.yang@mail.xjtu.edu.cn}
\\\vspace{0.2cm}
$^b$IBISC, Univ \'{E}vry, University of Paris-Saclay, \'{E}vry 91025, France\\
{\footnotesize E-mail address: feng.chu@univ-evry.fr}
\\\vspace{0.2cm}
$^c$State Key Lab for Manufacturing Systems Engineering, Xi'an 710049, China
\\\vspace{0.2cm}
$^d$Centre Excellence of Supply Chain, KEDGE Business School, 33400 Talence, France\\
{\footnotesize E-mail addresses: junkai.he@kedgebs.com}
}
\end{center}

\date{} 

\begin{abstract}
 In the context of fault-detection problems, the objective is to identify all defective items among a set of $n$ binary-state items using the minimum number of tests. The {group testing} paradigm, which allows testing a subset of items in a single test, serves as a fundamental technique for efficiently classifying large populations.
 We study a central problem in the combinatorial group testing model where the number $d$ of defective items is unknown in advance.
 Let $M_\alpha(d|n)$ denote the maximum number of tests required by an algorithm $\alpha$ for this problem, and $M(d,n)$ denote the minimum number of tests required in the worst case when $d$ is known in advance. An algorithm $\alpha$ is called a $c$-\emph{competitive algorithm} if there exist constants $c$ and $a$ such that, for $0\le d < n$, $M_{\alpha}(d|n)\le cM(d,n)+a$. 
 We design a new adaptive algorithm with a competitive constant $c \le 1.431$, thus pushing the competitive ratio below the best-known one of $1.452$. To achieve this, we propose a novel solution framework based on an unexplored up-zig-zag strategy and a studied strongly competitive algorithm.
  
  \

\end{abstract}

\section{Introduction} \label{sec:introduction}

As for many scenarios, the imperative task of screening a population of items to classify each item as either \emph{defective} or \emph{good} based on a binary characteristic holds paramount importance \citep{aprahamian2019optimal}. Intuitively, it is always too costly for large populations and limited resources for individual testing. Accordingly, \cite{dorfman1943detection} proposed a sequential testing strategy, called \emph{group testing}, to screen military inductees for syphilis in an economical manner. Roughly speaking, this scheme includes two stages: (\rmnum{1}) Test items in groups. The group is said to be \emph{contaminated} if it contains at least one defective item; otherwise, the group is said to be \emph{pure}. If a group is tested pure, then testing stops, and all items in that group are classified as good; otherwise, testing proceeds to the second stage. (\rmnum{2}) Test all items in the group individually and classify them according to their individual test outcome. Naturally, the choice of group size (i.e., the number of items in each group) has a significant impact on the efficiency of group testing. Group testing has been extensively studied but mostly under restrictive assumptions (as considered in this paper), such as \emph{perfect} tests (i.e., no classification errors), \emph{homogeneous} items (i.e., the probability of having the binary characteristic is the same for all items), and an infinite testing population, which leads to a focus on minimizing the number of tests.

Due to its simplicity and effectiveness, the problem later found various applications, including quality control in industrial product testing \citep{sobel1959group}, sequential screening of experimental variables \citep{li1962sequential}, multiaccess communication networks \citep{berger1984random, wolf1985born, goodrich2008improved}, DNA library screening \citep{pevzner1994towards}, blood screening for HIV tests \citep{wein1996pooled, zenios1998pooled}, data compression \citep{hong2002group}, testing and diagnosis for digital logic systems \citep{kahng2006new}, fault diagnosis in optical networks \citep{harvey2007non}, network security \citep{thai2011group}, failure detection in wireless sensor networks \citep{lo2013efficient}, etc. It is worth noting that among its recent applications, group testing has been widely advocated for expediting mass testing for COVID-19 PCR-based tests on a global scale \citep{ben2020large}. For an extensive summary of the work on this subject before 2000, the reader is referred to the monograph by \cite{du2000combinatorial}. In addition to the aforementioned applications, many variations of the classical group testing model have also been proposed to best adapt to the characteristics of different practical applications. The following list gives some examples that are not meant to be exhaustive: \citeauthor{de2003constructions} (\citeyear{de2003constructions}, \citeyear{de2006optimal}), \cite{bar2007applications}, \cite{claeys2010queueing}, \cite{feng2010efficient}, \cite{damaschke2013two}, \cite{kagan2014group}, \cite{zhang2015turnpike}, \cite{de2016constraining}, \cite{aprahamian2020optimal}, \citeauthor{cheng2019class} (\citeyear{cheng2019class}, \citeyear{cheng2021optimal}).

Let $M(d,n)$ denote the minimum number of tests required in the worst case if the number of defectives $d$ is known in advance. There are many studies on the bounds of $M(d,n)$ in the literature. \cite{hwang1972method} proposed a generalized binary-splitting algorithm that uses at most $\lceil \log\binom{n}{d} \rceil + d - 1$ tests in the worst case. Across the paper, $\log$ is of base 2 if no base is specified. As to lower bounds on $M(d,n)$, in general, the information-theoretic bound $\lceil \log\binom{n}{d} \rceil$ is still currently best known. The inequality is proved as Lemma 4.1.1 in \citeauthor{du2000combinatorial} (\citeyear{du2000combinatorial}, Chap. 4): For $0<\rho <1$, if $0< d <\rho n$, then
\begin{eqnarray}
  \log\binom{n}{d} \ge d\Big(\log \frac{n}{d} + \log(e\sqrt{1-\rho}) \Big) - 0.5\log d - 0.5\log(1-\rho) - 1.567.
\end{eqnarray}

In many applications, we probably only know the existence of defectives without any further information. In such situations, an algorithm whose number of tests used in the worst case will not be too much more than $M(d,n)$ is still desirable, regardless of the number of defectives. Motivated by the study of online problems \citep{MANASSE1990208, sleator1985amortized}, \cite{du1993competitive} proposed the concept of competitive algorithms for the group testing problem, then followed by \cite{bar1994new}, \cite{1994Modifications}, \cite{Schlaghoff2005Improved}, \cite{2014A}, \cite{wu2022improved}, and \cite{scheidweiler2023upper}. All the discussions that follow, as well as the main contribution of this paper, are presented in the context of the \emph{combinatorial group testing} (abbreviated as CGT) model. Beyond the classical CGT framework, several important extensions have been widely explored in the literature. For instance, as for the \emph{probabilistic group testing} model in which each item in $S$ is defective with some probability $p$, for the case of small defect probability $p\ll 1$ and large number $n$ of items, \cite{mezard2011group} obtained asymptotically optimal two-stage testing strategies. In addition, noisy group testing and testing under various constrained settings have been extensively studied \citep{aldridge2019group}.

Next, we follow two notions of competitiveness commonly used in the literature for the CGT with an unknown number $d$ of defectives (abbreviated as CGT-U). Let $A$ be an algorithm solving the CGT-U, and let $M_A(d|n)$ denote the number of tests performed by $A$ in the worst case. Then $A$ is called a $c$-competitive algorithm if there exists a constant $a$ such that for $0\le d < n$, $M_A(d|n)\le c M(d,n)+a$. Since $M(n,n) = 0$, and $M_A(n|n)\ge n$ for any $A$, the case $d=n$ is excluded in the definition. If $c$ is a constant, a $c$-competitive Algorithm $A$ is simply called a \emph{competitive} algorithm, and $c$ is called the \emph{competitive ratio} of $A$. \cite{du1993competitive} developed a $2.75$-competitive algorithm. The competitive ratio was improved to $2$ by \cite{bar1994new} and then to $1.65$ by \cite{1994Modifications}. \cite{Schlaghoff2005Improved} proposed an Algorithm $A_\epsilon$ satisfying $M_{A_\epsilon} \le (1.5+\epsilon)M(d,n)+\lceil{0.4/\epsilon}\rceil$ for all positive $\epsilon < 0.01$. Recently, \cite{scheidweiler2023upper} developed an algorithm with $c < 1.452$, which pushes the competitive ratio down under $1.5$.

\cite{du1994competitive} introduced the other notion of competitiveness. For an Algorithm $A$ solving the CGT-U, define
\begin{eqnarray*}
  n(d,k) = \max\{n|M(d,n) \le k\},~{\rm and}~ n_A(d|k) = \max\{n|M_A(d|n) \le k\}.
\end{eqnarray*}
Algorithm $A$ is called a \emph{strongly c-competitive} algorithm if there exists a constant $a$ such that for every $d\ge 1$ and $k\ge 1$, $n(d,k)\le c n_A(d|k) + a$. If $c$ is a constant, then a strongly $c$-competitive algorithm is also called a \emph{strongly competitive} algorithm. \cite{du1994competitive} proved that every strongly competitive Algorithm $A$ is also competitive in the original sense, as defined in \cite{du1993competitive}, and satisfies the following desirable property: $\lim_{n\rightarrow \infty}(M_A(d|n)/M(d,n))=1$. Moreover, they proposed an Algorithm \emph{DP} with $M_{DP}(d|n)\le d\log(n/d)+4d$ for $1\le d\le n$. \cite{Schlaghoff2005Improved} proved that there exists an Algorithm $B$ satisfying for $1\le d\le n$, $M_B(d|n) \le \lceil \log\binom{n}{d} \rceil + 2d$; however, the Algorithm $B$ implies by their constructive proof is rather complicated. \cite{2014A} presented a deterministic, strongly competitive Algorithm \emph{Z} based on a zig-zag approach with $M_Z(d|n)\le d\log(n/d)+3d+O(\log^2d)$. Recently, by using an innovative approach to handle contaminated subsets and to specify sizes of subsets to be tested, Algorithm $Z^*$ proposed by \cite{wu2022improved} improved the coefficient of term $d$ in the upper bound from $3$ to $5-\log 5$.

In this study, the focus lies on the competitive ratio. We introduce a novel solution framework for the CGT-U, which integrates our newly proposed \emph{up-zig-zag} strategy and a strongly competitive algorithm \citep{wu2022improved}. The solution framework comprises the following two main steps: (\rmnum{1}) conducting a small number of preliminary tests to determine whether there exists at least $x$ defective item(s) in the entire set, and (\rmnum{2}) selecting an appropriate procedure, either the up-zig-zag strategy or the strongly competitive algorithm, for the remaining items, depending on the outcome of the preliminary tests. Building upon this framework, we develop a new deterministic $c$-competitive Algorithm $Z^c$ and establish the following theorem:
\begin{theorem}
    \label{theo:paper}
    \centering
    For $0\le d < n, M_{Z^c}(d|n) \le 1.431 M(d,n) + 39$.
\end{theorem}

Rather than focusing solely on a single competitiveness perspective in previous studies, we conduct a comprehensive investigation into the effective integration of two methodological paradigms (i.e., $c$-competitive and strongly competitive), thereby leading to the above novel solution framework. It is worth noting that the up-zig-zag strategy itself can achieve a competitive ratio of $1.431$ for the CGT-U when $d = 0$ or $3 \le d < n$. Moreover, this strategy is straightforward to implement and performs well in practice. Figure~\ref{fig:Literature} provides a comparison between the existing results in the literature and our contributions to the CGT-U.
 
\begin{figure}[!ht]
  \centering
\begin{tikzpicture}
  \draw[{-Stealth[]}, dashed, line width=1pt] (-0.6,0) -- (14,0) node[below] {\tiny \emph{strongly competitiveness}};
  \draw[{-Stealth[]}, line width=1pt] (-0.6,0.2) -- (14,0.2) node[above] {\tiny \emph{c- competitiveness}};
  \draw[] (-0.6,-0.075) -- (-0.6,-0.075) node[above] {\Large {...}};
  \draw[] (0,-0.12) -- (0,-0.12) node[above] {\tiny {1993}};
  \draw[] (0.58,-0.12) -- (0.58,-0.12) node[above] {\tiny {1994}};
  \draw[] (4.1,-0.12) -- (4.1,-0.12) node[above] {\tiny {2005}};
  \draw[] (7.3,-0.12) -- (7.3,-0.12) node[above] {\tiny {2014}};
  \draw[] (9.9,-0.12) -- (9.9,-0.12) node[above] {\tiny {2022}};
  \draw[<-{Circle[]}, line width=1pt] (0,0.2) -- (0,4) node[right] {\scriptsize $2.75M(d,n)+O(1)$};
  \draw[] (0,0.2) -- (0,3.6) node[right] {\scriptsize (\citeauthor*{du1993competitive}, \citeyear{du1993competitive})};
  \draw[<-{Circle[]}, line width=1pt] (0.33,0.2) -- (0.33,2.8) node[right] {\scriptsize $2M(d,n)+5$};
  \draw[] (0.33,0.2) -- (0.33,2.4) node[right] {\scriptsize (\citeauthor{bar1994new},\citeyear{bar1994new})};
  \draw[<-{Circle[]}, line width=1pt] (0.33,0.2) -- (0.33,1.8) node[right] {\scriptsize $1.65M(d,n)+10$};
  \draw[-] (0.33,0.2) -- (0.33,1.4) node[right] {\scriptsize (\citeauthor{1994Modifications}, \citeyear{1994Modifications})};
  \draw[<-{Circle[]}, line width=1pt] (4.1,0.2) -- (4.1,1.6) node[right] {\scriptsize $(1.5+\epsilon)M(d,n)+\lceil 0.4/\epsilon\rceil$ for $\epsilon\in (0,0.01)$};
  \draw[-] (4.1,0.2) -- (4.1,1.2) node[right] {\scriptsize (\citeauthor*{Schlaghoff2005Improved}, \citeyear{Schlaghoff2005Improved})};
  \draw[<-{Circle[]}, line width=1pt] (10.05,0.2) -- (10.05,1.5) node[right] {};
  \draw[ ] (10.05,1.9) node[right] {\scriptsize {$1.452M(d,n)+O(1)$}};
  \draw[-] (10.05,0.2) -- (10.05,1.5) node[right] {\scriptsize (\citeauthor*{scheidweiler2023upper}, \citeyear{scheidweiler2023upper})};
  \draw[<-{Circle[]}, line width=1pt] (10.05,0.2) -- (10.05,1.1) node[right] {\scriptsize\bf {$1.431M(d,n)+39$}};
  \draw[-] (10.05,0.2) -- (10.05,0.7) node[right] {\scriptsize ({\bf this paper})};
  \foreach \x in {}
      \draw (\x,-0.095) -- (\x,-0.095) node[below] {\x};
  \draw[dashed] (0.33,0) -- (0.33,-2) node[right] {\scriptsize $d\log(n/d)+4d$};
  \draw[<-{Circle[open]}, dashed, line width=1pt] (0.33,0) -- (0.33,-2.4) node[right] {\scriptsize (\citeauthor*{du1994competitive}, \citeyear{du1994competitive})};
  \draw[dashed, line width=1pt] (4.1,0) -- (4.1,-1.8) node[right] {\scriptsize $\lceil d\log(n/d)\rceil+2d$};
  \draw[<-{Circle[open]}, dashed, line width=1pt] (4.1,0) -- (4.1,-2.2) node[right] {\scriptsize (\citeauthor*{Schlaghoff2005Improved}, \citeyear{Schlaghoff2005Improved})};
  \draw[dashed, line width=1pt] (7.3,0) -- (7.3,-1.4) node[right] {\scriptsize $d\log(n/d)+3d+O(\log^2d)$};
  \draw[<-{Circle[open]}, dashed, line width=1pt] (7.3,0) -- (7.3,-1.8) node[right] {\scriptsize (\citeauthor{2014A},\citeyear{2014A})};
  \draw[dashed, line width=1pt] (9.9,0) -- (9.9,-0.6) node[right] {\scriptsize $d\log(n/d)+(5-\log5)d+O(\log^2d)$};
  \draw[<-{Circle[open]}, dashed, line width=1pt] (9.9,0) -- (9.9,-1) node[right] {\scriptsize (\citeauthor{wu2022improved}, \citeyear{wu2022improved})};
\end{tikzpicture} 
\caption{.~{\bf A Summary of Known Results for the Combinatorial Group Testing Model on $c$-Competitive and Strongly Competitive Algorithms}\\ 
\emph{Note.}~Between the two $x$-axes, the upper one (solid line) corresponds to $c$-competitive algorithms, while the lower one (dotted line) corresponds to strongly competitive algorithms.}
\label{fig:Literature}
\end{figure}

\

\emph{Organization}. Recall some previous algorithms and results in \S~\ref{subsec:Preliminaries}, which will be used for subsequent algorithm development and analysis. Then detail the up-zig-zag strategy and the symmetric combinatorial algorithm in \S~\ref{sec:algorithm}. We go on to analyze the competitive ratio of our new algorithm in \S~\ref{sec:CR}. To maintain the coherence of the competitive ratio analysis, we place the proofs of two crucial lemmas (i.e., the parts that detail how to compute the two upper bounds on the number of total tests performed by strategy up-zig-zag) in \S~\ref{sec:upzigzag}. Finally, conclude in \S~\ref{sec:Cr}. For the sake of clarity and continuity, some proofs of claims and lemmas are relegated to Appendix A.

\section{Preliminaries}
\label{subsec:Preliminaries}

For the sake of convenience, let a multiset $(S_n, d|n)$ be a general CGT-U instance consisting of  $n$ items, where the number $d$ ($0 \le d \le n$) of defective items is unknown in advance:
\begin{eqnarray*}
(S_n, d|n)=\big\{\underbrace{1,1,\cdots,1}_d,\underbrace{0,0,\cdots,0}_{n-d}\big\}\Big|^{\rm 0:~good~item}_{\rm 1:~defective~item}
\end{eqnarray*} 
To further develop the foundation for this work, we next introduce some key components that will be used for designing our newly proposed algorithm.

A binary splitting procedure DIG (see Procedure~\ref{alg:DIG}) to identify one defective (and an unspecified number of good items) from a contaminated set $X$ is described first. We use $\lceil x\rceil$ to denote the ceiling function of $x$.

\begin{algorithm}
  \renewcommand{\algorithmcfname}{Procedure}
  \small{
  \caption{DIG$(X)$}
  \label{alg:DIG}
  \KwIn{a contaminated set $X$;}
  \KwOut{one identified defective and an unspecified number of good items in $X$;}
  \eIf{$|X|>1$}
  {
    \rm Test $\{X'|~X'\subseteq X, |X'|\leftarrow \lceil |X|/2\rceil\}$\;
      \eIf{$X'$ is contaminated}
      {
        \rm $X\leftarrow X'$\;
      }
      {
        \rm Identify all items in $X'$ as good items\;
        \rm $X\leftarrow X\setminus 
        X'$\;
      }
  }
  {
    \rm Identify the single item in $X$ as defective.
  }
  }
\end{algorithm}

The following folklore lemma about procedure DIG will be used for later analysis.

\begin{lemma}\label{lem:DIG}
For a contaminated set $X$ with $2^{i-1}<|X|\le2^i$, where $i\ge1$ is an integer, the total number of tests performed by procedure DIG on $X$ is $i$.
\end{lemma}

We then describe an integer sequence, a procedure $4$-Split, and an improved zig-zag approach (here denoted by $Z^d$), which were proposed by \cite{wu2022improved}. The integer sequence is defined as follows: $a_0=1, a_1=2, a_2=3, \cdots, a_i=\lceil 3\cdot 2^{i-2}\rceil$ for $i \ge 0$. Procedure $4$-Split (see Procedure~\ref{alg:4Split}) will be applied for later procedures to identify one defective item (and an unspecified number of good items) from a contaminated set $X$ of size at most $a_k$ items. Intuitively, the core idea of this procedure is to partition a contaminated set into four disjoint subsets of a certain size for testing.

\begin{algorithm}
  \renewcommand{\algorithmcfname}{Procedure}
\small{
  \caption{4-Split$(X,k)$ \citep{wu2022improved}}
  \label{alg:4Split}
  \KwIn{a contaminated set $X$, and an integer $k\ge 0$ such that $|X|\le a_k$;}
  \KwOut{one identified defective item and all good items in $X$ identified before it;}
  \If{$|X|=1$}
  {
    \rm Identify the single item in $X$ as defective without any test, remove this item from $X$\; \label{S'=1}
  }
  \ElseIf{\rm $2\le |X|\le 3$}
  {
    \rm Perform individual tests on items in $X$ until the first defective item is found, remove this defective item (and all good items identified before it) from $X$\; \label{S'<=3}
  }
  \Else
  {
    \rm Partition $X$ into four pairwise disjoint sets $Y,~Z,~U,~V$. such that
    \begin{equation*}
    \begin{aligned}
        |Y|&=\min \{2^{k-2}, |X|\},\\
        |Z|&=\min \{2^{k-2}, |X|-|Y|\},\\
        |U|&=\min \{2^{k-3}, |X|-|Y|-|Z|\},\\
        |V|&=\min \{2^{k-3}, |X|-|Y|-|Z|-|U|\}\\
    \end{aligned}
    \end{equation*}
    \eIf{$Z,U$, and $V$ are all empty}
    {
      Apply {\bf DIG}($Y$) to identify one defective item in $Y$, and remove this item from $X$\;
    }
    {
      Test $Y, Z, U$, and $V$ sequentially until the first contaminated subset $X'$ is found. Apply {\bf DIG}($X'$) to identify one defective item in $X'$. Remove this defective item (and all good items tested before $X'$, if any) from $X$\;
    }
  } 
  }
\end{algorithm}

Procedure~$Z^d$ is an essential part of our new algorithm. We present its framework in Procedure~\ref{alg:Zd}. Each time, the algorithm chooses a subset $S'$ ($|S'|\le a_k$; if the size of the entire set is less than $a_k$, then select it as $S'$), to identify a defective item from it or to determine that the subset is pure. The size of the current subset depends on the size and test result of the preceding subset. Roughly speaking, each time the size is increased if the preceding subset is tested to be pure, and the size is reduced if the preceding subset is tested to be contaminated. The process continues until every item in the entire set is successfully identified.

\begin{algorithm}
  \renewcommand{\algorithmcfname}{Procedure}
  \small{
  \caption{Improved Zig-Zag Approach (termed $Z^d$) \citep{wu2022improved}}
  \label{alg:Zd}
  \KwIn{an instance $(S_n, d|n)$, and an integer $k = \lceil\log (4n/3)\rceil$;}
  \KwOut{a complete identification of all defective and good items in $S_n$;}
  \While{$S_n\neq \emptyset$}
  {
    \rm Test \big\{\rm $S'\big|~S'\subseteq S_n, |S'|\leftarrow \min\{a_k, |S_n|\}$\big\}\;
    \eIf{$S'$ is pure}
    {
      \rm Identify all items in $S'$ as good items and remove them from $S_n$\; 
      \rm $k\leftarrow k + 1$\;
    }
    {
      \eIf{$k>0$}
      {
        \rm Apply {\bf 4-Split}$(S',k)$ to identify one defective item and a certain number of good items in $S'$, and remove these items from $S_n$\;
        \rm $k\leftarrow k - 1$\;
      }
      {
        \rm Identify the single item in $S'$ as defective and remove 
        it from $S_n$\;
        \rm $k\leftarrow k$\;
      }
    }
  }
  }
\end{algorithm}

The following lemma (Theorem 3.12 in \cite{wu2022improved}) about procedure $Z^d$ is stated here without giving proof.

\begin{lemma}\label{lem:Z1mian} \citep{wu2022improved} For an instance $(S_n,d|n)$, $M_{Z^d}(d|n)$, the total number of tests used by $Z^d$ in the worst case satisfies that
$$
M_{Z^d}(d|n)\le d\log \frac{n}{d}+ (5-\log 5)d + 0.5\log^2 d + \big(\log(5/3) + 1.5\big)\log d + 4.
$$
\end{lemma}

As a consequence of Lemma~\ref{lem:Z1mian}, the following two corollaries for later analysis of our new algorithm can be deduced.

\begin{corollary}\label{coro:Z1corollary1}
  Given an instance $(S_n, d|n)$, such that $S_n$ contains at least $3n/4$ good items. Let $\psi$ be some absolute positive constant. Suppose that at least $3n/4$ good items can be detected by at most $\psi$ tests first, and then apply procedure $Z^d$ to the remaining items. As a result, the total number of tests for $(S_n, d|n)$, $M(d|n)$, satisfies that
  $$
  M(d|n) \le 1.431
   d(\log\frac{n}{d} + 1) + 4 + \psi.
  $$
\end{corollary}

\noindent{\bf Proof.}\, By Lemma~\ref{lem:Z1mian}, we have
\begin{eqnarray*}
  M(d|n) &=& M_{Z^d}(d|n/4) + \psi\\
  &\le& d\log \frac{n/4}{d}+(5-\log 5)d+0.5\log^2 d+\big(\log(5/3) + 1.5\big)\log d+4+\psi\\
  &=& d\log \frac{n}{d}+(3-\log 5)d+0.5\log^2 d+\big(\log(5/3) + 1.5\big)\log d+4+\psi.
\end{eqnarray*}
Define $\Delta_f = 4 + \psi$ and let $f(d) = 1.431 d(\log\frac{n}{d} + 1) + \Delta_f - M(d|n)$, then we have
\begin{eqnarray}
  f(d) &\ge& 0.431d\log\frac{n}{d} + (\log 5-1.55)d - 0.5\log^2d - \big(\log(5/3)+1.5\big)\log d - 4 - \psi + \Delta_f\nonumber\\
  &\ge& (\log5-0.688)d - 0.5\log^2d - \big(\log(5/3)+1.5\big)\log d\label{Z1-0}\\
  &=& g(d),\nonumber
\end{eqnarray}
where inequality~(\ref{Z1-0}) holds for $d\le n/4$.

First, we get $g'(d) = \frac{(\log5-0.688)\cdot \ln 2\cdot d - \log d - \log(5/3)-1.5}{\ln 2\cdot d}$. The calculation shows that $g'(d)<0$ when $0\le d\le 3$ and $g'(d)>0$ when $d\ge 4$. Thus, 
$$
g(d)\ge \min\{g(3),g(4)\} > 0,
$$
which implies $f(d)\ge 0$. Therefore, $ M(d|n) \le 1.431d\big(\log\frac{n}{d}+ 1\big) + 4 + \psi$. The corollary is proven. \hfill $\Box$

\begin{corollary}\label{coro:Z1corollary2}
  Given an instance $(S_n, d|n)$, such that $S_n$ contains at least $3n/4$ good items. Let $\psi$ be some absolute positive constant. Suppose that at least $3n/4$ good items can be detected by at most $\psi$ tests first, and then apply procedure $Z^d$ to the remaining items. As a result, the total number of tests for $(S_n, d|n)$, $M(d|n)$, satisfies that
  $
  M(d|n) \le 1.07325n + 4 + \psi.
  $
\end{corollary}

\noindent{\bf Proof.}\, Based on Corollary~\ref{coro:Z1corollary1}, we let
$$
f(d) = 1.431n - M(d|n) \ge 1.431n - 1.431d(\log\frac{n}{d} + 1) - 4 - \psi = g(d)
$$

First, $g'(d) = 1.431[\log d-\log (2n/e)]$. When $d=2n/e$, $g'(d) = 0$. Since $1\le d \le n/4$, we have
\begin{eqnarray*}
  g(d) &\ge& g(n/4) \\
  &=& 1.431n - 1.431\cdot (n/4) \cdot (\log\frac{n}{n/4} + 1) - 4 - \psi \\
  &=& 0.35775n - 4 - \psi.
\end{eqnarray*}
Then, $1.431n-M(d|n) \ge 0.35775n - 4 - \psi$. Thus, $M(d|n) \le 1.07325n + 4 + \psi$. \hfill $\Box$

\

Finally, we need the following result proposed by \cite{1994Modifications}.

\begin{lemma}\label{lem:com-d}
\citep{1994Modifications}
  Let $d=d_1+d_2$ and $n=n_1+n_2$ where $d_1\ge0, d_2\ge0, n_1>0, n_2>0$. Then
  $$
  d_1\log\frac{n_1}{d_1}+d_2\log\frac{n_2}{d_2}\le d\log\frac{n}{d}.
  $$
\end{lemma}

\section{The Symmetric Combinatorial Algorithm: Overview and Intuition}\label{sec:algorithm}

Formally, our new adaptive algorithm (here denoted by $Z^{c}$ in short) is essentially a combination of two analogous strategies. To better introduce the new algorithm later, we first explain the simple intuition that lies behind it. In the beginning, we propose an up-zig-zag strategy (or $Z^u$ in short) that can solve the CGT-U. However, the limitation of this new strategy is that it achieves a competitive ratio of $1.431$ only when $d=0$ and $3 \le d < n$. To overcome this limitation, we couple procedure $Z^u$ with procedure $Z^d$ to form Algorithm~$Z^{c}$, thereby generalizing the competitive ratio of $1.431$ to the case of $0\le d < n$. Given the similarities in mechanism between procedures $Z^d$ and $Z^u$, we name Algorithm $Z^{c}$ a symmetric combinatorial algorithm. 

In the following, we present details of the up-zig-zag strategy and Algorithm $Z^{c}$.

\subsection{Framework of Strategy Up-Zig-Zag}\label{subsec:framework_upzigzag}

The basic framework of the up-zig-zag strategy is as follows. For an instance $(S_n,d|n)$, each time we choose a subset $S'$ of size $a_k$ to identify it as a pure set or a contaminated set. If $S'$ is pure, increase $k$ by one; otherwise, identify one defective item from it by procedure 4-Split and decrease $k$ by one (while $k=0$, $k$ is unchanged). In particular, we log the number of subsets that have been consecutively tested as pure sets and assign the value to $p_c$. If $p_c$ equals $6$ and the size of the remaining entire subset $S\subseteq S_n$ is larger than $a_{k}$, test $S$ as a group. If the test result of $S$ is pure, terminate; otherwise, we go on to test the next subset of size $a_{k+1}$. The above process continues until every item in $S_n$ is successfully identified. In the beginning, we test a subset $S'$ of size $a_k$ with $k = 0$.

Whereas after conducting preliminary calculations, we found a barrier that the number of tests performed on a contaminated set $S'$ of size $2$ or $3$ by the procedure 4-Split alone can not reach an expected upper bound. Thus, to address the challenge, we design two ad hoc refined schemes: procedures 2-Test and 3-Test (see Procedures~\ref{strategy:S2} and~\ref{strategy:S3}), for such $S'$. Procedure 2-Test is used to handle a set $S^2$ of size $2$ containing at least one defective item. First, test items in $S^2$ individually. And then, we introduce an auxiliary binary variable $\Delta_{\rm TS}$(= $0$ or $1$) to distinguish the different test outcomes on $S^2$. That is, only if $S^2$ contains two defective items, i.e., $S^2=\{1,1\}$, decrease $k$ by one and set $\Delta_{\rm TS}=0$; otherwise, increase $k$ by one and set $\Delta_{\rm TS}=1$. For a set $S^3$ of size $3$ containing at least one defective item, procedure 3-Test tests all items in it individually only if $\Delta_{\rm TS}$ is equal to $1$ before testing. Indeed, the above two refinement schemes are the key to deriving a desirable upper bound for $M_{Z^u}(d|n)$ that depends only on $n$, see Lemma~\ref{lem:UpperBound2-Zu}.

\begin{algorithm}
  \renewcommand{\algorithmcfname}{Procedure}
  \small{
  \caption{2-Test$(S',k,p_c,\Delta_{\rm TS})$}
  \label{strategy:S2}
  \KwIn{a contaminated set $S'$ of size $2$, and three integer variables $k = 1, p_c, \Delta_{\rm TS} = 0$;}
  \KwOut{updated values of $k$, $p_c$, and $\Delta_{\rm TS}$, and a complete identification of all defective and good items in $S'$;}
  \rm Test all items in $S'$ individually\;
  \eIf{$S'=\{0,1\}$}
    {
      \rm Identify one defective item and one good item in $S'$, and remove these items from $S_n$\;
      \rm $k\leftarrow 2$, $p_c\leftarrow p_c + 1$, and  $\Delta_{\rm TS} \leftarrow 1$\;
    }
    {
      \rm Identify two defective items in $S'$, and remove these items from $S_n$\;
      \rm $k\leftarrow 0$, $p_c\leftarrow 0$, and $\Delta_{\rm TS}\leftarrow 0$\;
    }
  }
\end{algorithm}

\begin{algorithm}
  \renewcommand{\algorithmcfname}{Procedure}
  \small{
  \caption{3-Test$(S',k,p_c,\Delta_{\rm TS})$}
  \label{strategy:S3}
  \KwIn{a contaminated set $S'$ of size $3$, and three integer variables $k=2, p_c, \Delta_{\rm TS}=1$;}
  \KwOut{updated values of $k$, $p_c$, and $\Delta_{\rm TS}$, and a complete identification of all defective and good items in $S'$;}
  \rm Test all items in $S'$ individually\; 
  \rm Identify a certain number of defective and good items from $S'$, and remove them from $S_n$\;
  \rm $k\leftarrow 1$, $p_c\leftarrow 0$, and $\Delta_{\rm TS}\leftarrow 0$\;
  }
\end{algorithm}

The basic framework of the up-zig-zag strategy, coupled with procedures~2-Test and 3-Test, can lead up to procedure $Z^u$ as depicted in Procedure~\ref{alg:Zu}. Intuitively, the process of procedure $Z^u$ shares similarities with procedure $Z^d$, changing the size of the current test subset according to the outcome of the preceding test. However, the latter approach selects the entire set $S_n$ as the first test subset. Conversely, procedure $Z^u$ starts with a subset of size $1$. From the perspective of the whole testing process, procedure $Z^d$ has a decreasing trend in the size of each test subset, but our method does the opposite, so we call procedure $Z^u$ the up-zig-zag strategy to distinguish it from the original zig-zag approach. Compared to procedure $Z^d$, procedure $Z^u$ also incorporates several additional specialized operations. These notable differences make the analysis method for $Z^d$ no longer applicable to $Z^u$.

\begin{algorithm}
  \renewcommand{\algorithmcfname}{Procedure}
  \small{
  \caption{Up-Zig-Zag Strategy ($Z^u$)}
  \label{alg:Zu}
  \KwIn{an instance $(S_n, d|n)$, a set $S'$, and three integer variables: $k$, $p_c$, $\Delta_{\rm TS}$;}
  \KwOut{a complete identification of all defective and good items in $S_n$;}
  Initially, set $S' = \emptyset$, $k = 0$, $p_c = 0$ and $\Delta_{\rm TS} = 0$\;
  \While{$S_n\neq \emptyset$}
  {
   \eIf{$S'=\emptyset$}
    {
     \eIf{$p_c=6$ $\wedge$ $|S_n| > a_k$}
     {
      \rm Test $S_n$\;
      \eIf{$S_n$ is pure}
      {
        \rm Identify all items $S_n$ as good items and remove them from $S_n$\;
      }
      {
        \rm $S'\leftarrow \max\big\{X|~X \subseteq S_n, |X|\le a_k \big\}$\;
      }
     }
      {
        \rm $S'\leftarrow \max\big\{X|~X \subseteq S_n, |X|\le a_k \big\}$\;
      }
    }
    {
      \rm Test $S'$\;
      \eIf{$S'$ is pure}
      {
        \rm Identify all items in $S'$ as good items and remove them from $S_n$\;
        \rm $k\leftarrow k + 1$ and $p_c\leftarrow p_c + 1$\;
      }
      {
        \rm Check the value of $k$\;
        \eIf{$k=1$}
        {
          \rm {\bf Procedure~2-Test}$(S',k,p_c,\Delta_{\rm TS})$\;
        }
        {
          \rm Check the value of $\Delta_{\rm TS}$\;
          \eIf{$\Delta_{\rm TS}=1$ $\wedge$ $k=2$}
          {
            \rm {\bf Procedure~3-Test}$(S',k,p_c,\Delta_{\rm TS})$\; 
          }
          {
            \rm $p_c\leftarrow 0$\;
          \eIf{$k=0$}
          {
            \rm Identify one defective item in $S'$ and remove it from $S_n$\;
            \rm $k\leftarrow k$\;
          }
          {
            \rm Apply {\bf 4-Split}$(S',k)$ to identify one defective item and a certain number of good items in $S'$, and remove these items from $S_n$\;
            \rm $k\leftarrow k-1$\;
          }
          }
        }
      }
    }
    $S'\leftarrow \emptyset$\;
  }
  }
\end{algorithm}

\subsection{Framework of Symmetric Combinatorial Algorithm}\label{subsec:framework_alg}

Algorithm~\ref{alg:Zc} presents an overview of the symmetric combinatorial Algorithm $Z^{c}$. For a given instance $(S_n,d|n)$, the process begins with partitioning the initial set $S_n$ into four item-disjoint subsets of equal size and an unspecific subset $R_1$ of size $n_{R_1}$($\le 3$). Let $S_1$ be the collection of the above four equal-size subsets, such as $S_1:= \big\{S_1^v\big|~|S_1^v|=n_1,~\forall v=1,2,3,4 \big\}$. Accordingly, $S_n = S_1\cup R_1$ and $S_1\cap R_1 = \emptyset$. First, test items in $R_1$ individually to identify them all with no more than three tests. Then test each subset in $S_1$ simultaneously and record the number of contaminated subsets as $\alpha^1$. 

If $\alpha^1=0$, terminate; if $\alpha^1=1$, apply procedure $Z^d$ to it; and if $\alpha^1\ge3$, merge all contaminated subsets as a new subset and apply procedure $Z^u$ to it. However, it's more complicated if $\alpha^1=2$. Combine the resulting two contaminated subsets (defined as $S_1^1$ and $S_1^2$) as a whole and then partition it into four equal-size item-disjoint subsets which form a collection $S_2$, such as $S_2:= \big\{S_2^v\big|~|S_2^v|=n_2,~\forall v=1,2,3,4 \big\}$, and an unspecific subset $R_2$ of size $n_{R_2}$($\le 2$). So we have $S_1^1\cup S_1^2 = S_2\cup R_2$ and $S_2\cap R_2 = \emptyset$. We first test items in $R_2$ individually to detect them all with no more than two tests. Then test each subset in $S_2$ simultaneously and record the number of contaminated subsets as $\alpha^2$. If $\alpha^2 \ge 3$, merge all contaminated subsets into a new subset and apply procedure $Z^u$ to it; otherwise, merge all contaminated subset(s) into a new subset and apply procedure $Z^d$ to it. The process continues until every defective item in $S_n$ is successfully identified.

\begin{algorithm}
  \small{
  \renewcommand{\thealgocf}{1}
  \caption{Symmetric Combinatorial Algorithm ($Z^{c}$)}
  \label{alg:Zc}
  \KwIn{an instance $(S_n, d|n)$, and four nonnegative integers $n_1,n_{R_1},n_2$, and $n_{R_2}$;}
  \KwOut{a complete identification of all defective and good items in $S_n$;}
  Initially, set $n_1 \leftarrow \lfloor n/4\rfloor$, $n_{R_1} \leftarrow n - 4 n_1$. Partition $S_n$ into 4 pairwise disjoint subsets such that $\big\{S_1^v\big|~|S_1^v|=n_1 \big\}$ for $v\in \{1,2,3,4\}$, and a subset $\big\{R_1\big|~|R_1|=n_{R_1}\big\}$\;
  \While{$R_1\neq\emptyset$}
  {
    \rm Test all items in $R_1$ individually\;
    \rm Identify a certain number of contaminated items and pure items in $R_1$, and remove these items from $S_n$\;
  }
  \While{$\forall S_1^v \neq \emptyset$}
  {
    \rm Test $S_1^1,S_1^2,S_1^3$, and $S_1^4$ simultaneously\;
    \rm Record the number of contaminated subsets and assign the value to $\alpha^1$\;
    \If{$\alpha^1=0$}
    {
      \rm Identify all items as good items and remove them from $S_n$;
    }
    \ElseIf{$\alpha^1=1$}
    {
      \rm Identify all items in the pure subset as good items and remove them from $S_n$\;
      \rm Apply {\bf Procedure}~$Z^d$ to $S_n$\;
    }
    \ElseIf{$\alpha^1=2$}
    {
      \rm Identify all items in each pure subset as good items and remove them from $S_n$\;
      Set $n_2 \leftarrow \lfloor n_1/2\rfloor$, $n_{R_2} \leftarrow 2n_1 - 4 n_2$\;
      Partition $S_n$ into 4 pairwise disjoint subsets such that $\big\{S_1^v\big|~|S_2^v|=n_2 \big\}$ for $v\in \{1,2,3,4\}$, and a subset $\big\{R_2\big|~|R_2|=n_{R_2}\big\}$\;
      \While{$R_2\neq \emptyset$}
      {
        \rm Test all items in $R_2$ individually\;
        \rm Identify a certain number of contaminated items and pure items in $R_2$, and remove these items from $S_n$\; 
      }
      \While{$\forall S_v^2\neq \emptyset$}
      {
        \rm Test $S_2^1,S_2^2$, $S_2^3$ and $S_2^4$ simultaneously\;
        \rm Record the number of contaminated subsets and assign the value to $\alpha^2$\;
        \eIf{$0<\alpha^2\le2$}
        {
          \rm Identify all items in each pure subset as good items and remove them from $S_n$\;
          \rm Apply {\bf Procedure}~$Z^d$ to $S_n$\;
        }
        {
          \rm Identify all items in each pure subset as good items and remove them from $S_n$\;
          \rm Apply {\bf Procedure}~$Z^u$ to $S_n$\;
        }
      }
    }
    \Else
    {
      \rm Identify all items in the pure subset as good items and remove them from $S_n$.\;
      \rm Apply {\bf Procedure}~$Z^u$ to $S_n$\;
    }
  }
  }
\end{algorithm}

\section{Competitive Ratio}\label{sec:CR}

In this section, we present the results on the competitive ratio of Algorithm~$Z^c$, which is bounded by $1.431$, as stated in Theorem~\ref{theo:paper}. To establish this result, we first derive two upper bounds on the number of tests performed by Procedure~$Z^u$, as stated in Lemma~\ref{lem:UpperBound1-Zu} and Lemma~\ref{lem:UpperBound2-Zu}. These bounds serve as the basis for establishing Theorem~\ref{theo:paper}. Since the proofs of these lemmas are rather involved, they are deferred to Section~\ref{sec:upzigzag} for clarity and continuity of exposition. 

\begin{lemma}
  \label{lem:UpperBound1-Zu}
  For an instance $(S_n, d|n)$, where $ 3 \le d \le n$, the total number of tests used by procedure $Z^u$ in the worst case satisfies that
  $$
  M_{Z^u}(d|n) \le 1.431 d \big(\log\frac{n}{d}+1.1242\big)+23.
  $$
\end{lemma}

The following lemma has the same premises as Lemma~\ref{lem:UpperBound1-Zu}, but a different upper bound on $M_{Z^u}(d|n)$ depending on $n$ alone.

\begin{lemma}
  \label{lem:UpperBound2-Zu}
  For an instance $(S_n, d|n)$, where $0\le d \le n$, the total number of tests used by procedure $Z^u$ in the worst case satisfies that $M_{Z^u}(d|n) \le 1.4n$.
\end{lemma}

Based on the above results, two upper bounds on $M_{Z^{c}}(d|n)$ can be derived. For brevity, all of the steps we referred to in Section~\ref{sec:CR} are in Algorithm~\ref{alg:Zc}. For the sake of calculation, we partition the entire process of Algorithm $Z^{c}$ into two parts: Part~\Rmnum{1}, including the process from Step $(2)$ to Step $(4)$; Part~\Rmnum{2}, including the process from Step $(5)$ to Step $(31)$. Let $t_{\rm\Rmnum{1}}$ and  $t_{\rm\Rmnum{2}}$ be the number of tests incurred by Part \Rmnum{1} and Part \Rmnum{2}, respectively. Thus, we conclude that the total number of tests used by Algorithm~$Z^{c}$, defined as $t_{Z^{c}}$, is the sum of $t_{\rm\Rmnum{1}}$ and $t_{\rm\Rmnum{2}}$, such as,
\begin{eqnarray}
  t_{Z^{c}} = t_{\rm\Rmnum{1}} + t_{\rm\Rmnum{2}}.\label{totaltestZD}
\end{eqnarray}

The following analysis is based on formula~(\ref{totaltestZD}), Lemma~\ref{lem:UpperBound1-Zu} and Lemma~\ref{lem:UpperBound2-Zu}.

\begin{lemma}
  \label{lem:UpperBound-Zd}
  For an instance $(S_n, d|n)$, Algorithm $Z^{c}$ satisfies: {\rm(\rmnum{1})}~$M_{Z^{c}}(d|n) \le 1.431 d\big(\log \frac{n}{d} + 1.1242\big) + 23$, and {\rm(\rmnum{2})}~$M_{Z^{c}}(d|n) \le 1.4n + 13$, for $0\le d \le n$.
\end{lemma}

{\noindent\bf Proof.}\, We first estimate two terms $t_{\rm \Rmnum{1}}$ and $t_{\rm \Rmnum{2}}$ separately, and then combine them to reach two upper bounds on $M_{Z^{c}}(d|n)$. 

Indeed, Part~\Rmnum{1} is the process of testing $R_1$. Let $n_{\rm \Rmnum{1}}$ and $d_{\rm \Rmnum{1}}$ be the total number of items and defective items in Part \Rmnum{1}, respectively. Consider two cases on $d_{\rm\Rmnum{1}}$: (\rmnum{1})~for $d_{\rm\Rmnum{1}}=0$, $t_{\rm\Rmnum{1}}\le n_{\rm\Rmnum{1}}$; (\rmnum{2})~for $d_{\rm\Rmnum{1}}\ge 1$, $t_{\rm\Rmnum{1}} \le 1.431d_{\rm\Rmnum{1}}\big(\log (n_{\rm\Rmnum{1}}/d_{\rm\Rmnum{1}}) + 1.1242\big) + n_{\rm\Rmnum{1}} - d_{\rm\Rmnum{1}}$. Since $n_{\rm \Rmnum{1}} = n_{R_1}\le 3$, the following two results on Part~\Rmnum{1} can be derived,
\begin{eqnarray}
  t_{\rm\Rmnum{1}} &\le& 1.431d_{\rm\Rmnum{1}}\big(\log \frac{n_{\rm\Rmnum{1}}}{d_{\rm\Rmnum{1}}} + 1.1242\big) + 3,\label{Part1}\\
  {\rm and~}t_{\rm\Rmnum{1}} &=& n_{\rm \Rmnum{1}} \le 1.4n_{\rm \Rmnum{1}}.\label{Part1-1}
\end{eqnarray}

Part~\Rmnum{2} is the process of testing all subsets in $S_1$. Let $n_{\rm \Rmnum{2}}$ and $d_{\rm \Rmnum{2}}$ be the total number of items and defective items in Part \Rmnum{2}, respectively. By Step (6) and Step (7), it first tests all subsets in $S_1$ simultaneously with four tests to identify the number of contaminated subsets and record the value as $\alpha^1$. By Algorithm~$Z^{c}$, different test strategies are applied to the remaining items according to different values of $\alpha^1$.

If $\alpha^1 = 0$, all items in Part~\Rmnum{2} are good and so $t_{\rm \Rmnum{2}} = 4$. Then we have
\begin{eqnarray}
  t_{\rm \Rmnum{2}} &\le& 1.431d_{\rm \Rmnum{2}}\big(\log \frac{n_{\rm \Rmnum{2}}}{d_{\rm \Rmnum{2}}} + 1.1242\big) + 4,\label{Part2-0}\\
  {\rm and~}t_{\rm\Rmnum{2}} &\le& 1.4n_{\rm \Rmnum{2}}.\label{Part2-0-1}
\end{eqnarray}

If $\alpha^1 = 1$, there exists only one contaminated subset in $S_1$, and so apply procedure $Z^d$ to it. From Step (11) to Step (12), it can be seen that the total number of items incurred by procedure $Z^d$ is $n_{\rm \Rmnum{2}}/4$. Testing out three-fourths of the items from $S_1$ as good by four tests implies that Corollary~\ref{coro:Z1corollary1} and Corollary~\ref{coro:Z1corollary2} apply in this case and $\psi=4$, so we have
\begin{eqnarray}
  t_{\rm \Rmnum{2}} &\le& 1.431d_{\rm \Rmnum{2}}\big(\log \frac{n_{\rm \Rmnum{2}}}{d_{\rm \Rmnum{2}}} + 1\big) + 4 + \psi ~<~ 1.431d_{\rm \Rmnum{2}}\big(\log \frac{n_{\rm \Rmnum{2}}}{d_{\rm \Rmnum{2}}} + 1.1242 \big) + 8,\label{Part2-1}\\
  {\rm and~}t_{\rm\Rmnum{2}} &\le& 1.0725n_{\rm \Rmnum{2}} + 4 + \psi ~<~ 1.4n_{\rm \Rmnum{2}} + 7,\label{Part2-1-1}
\end{eqnarray}
where (\ref{Part2-1}) holds by Corollary~\ref{coro:Z1corollary1} and (\ref{Part2-1-1}) holds by Corollary~\ref{coro:Z1corollary2}.

If $\alpha^1 = 2$, there exist two contaminated subsets in $S_1$. Partition these two contaminated subsets into four subsets which form $S_2$ and a remaining subset $R_2$, such that $2n_1=4n_2+n_{R_2}$. First, perform no more than two tests for $R_2$. Let $d_{R_2}$ be the number of defective items in $R_2$. Then test $S_2$ by procedure $Z^d$ or $Z^u$. Similar to $S_1$, it tests all subsets in $S_2$ simultaneously by four tests to identify the number of contaminated subsets and record the value as $\alpha^2$ as shown in Step (21) and Step (22). Apply two different group testing strategies to the remaining items based on different values of $\alpha^2$, considering the following Cases 1 and 2.

{{\bf Case 1.}}~$1\le \alpha^2\le 2$. If $\alpha^2 = 2$, we first merge the two existing contaminated subsets into a new set and then apply procedure $Z^d$ to it; otherwise, we apply procedure $Z^d$ to only one contaminated subset directly. Intuitively, $2n_1+2n_2$ good items can be tested out by only eight tests (four tests for $S_1$ and four tests for $S_2$) first, and the total number of items incurred by procedure $Z^d$ is no more than $2n_2$. Since
  \begin{eqnarray*}
    2n_2 &=& n_1 - \frac{n_{R_2}}{2} ~=~ \frac{n_{\rm \Rmnum{2}}}{4} - \frac{n_{R_2}}{2} ~\le~ \frac{1}{4}(n_{\rm \Rmnum{2}} - n_{R_2}), \\
    {\rm and~}2n_1+2n_2 &=& \frac{n_{\rm \Rmnum{2}}}{2} + \frac{n_{\rm \Rmnum{2}}}{4} - \frac{n_{R_2}}{2} ~\ge~ \frac{3}{4}(n_{\rm \Rmnum{2}} - n_{R_2}),
  \end{eqnarray*}  
Corollary~\ref{coro:Z1corollary1} and Corollary~\ref{coro:Z1corollary2} can apply in this case with $\psi = 8$. Thus, the total number of tests in Part~{\Rmnum{2}} satisfies that
  \begin{eqnarray*}
    t_{\rm \Rmnum{2}}
     \le 1.431(d_{\rm \Rmnum{2}}-d_{R_2})\big(\log \frac{n_{\rm \Rmnum{2}}-n_{R_2}}{d_{\rm \Rmnum{2}}-d_{R_2}} + 1\big) + 4 + \psi + 2 \le 1.431d_{\rm \Rmnum{2}}\big(\log \frac{n_{\rm \Rmnum{2}}}{d_{\rm \Rmnum{2}}} + 1.1242\big) + 14,
  \end{eqnarray*}
which holds by Corollary~\ref{coro:Z1corollary1} and Lemma~\ref{lem:com-d}, and
  \begin{eqnarray*} 
    t_{\rm\Rmnum{2}} &\le& 1.07325{n_{\rm \Rmnum{2}}} + 4 + \psi + 2 ~<~ 1.4n_{\rm \Rmnum{2}} + 13,
  \end{eqnarray*}
which holds by Corollary~\ref{coro:Z1corollary2} and the fact of $n_{\rm\Rmnum{2}}\ge 8n_2 + n_{R_2} \ge 8$.

{{\bf Case 2.}~$3 \le \alpha^2\le 4$.} So there exist at least three contaminated subsets in $S_2$, implying $d_{\rm \Rmnum{2}} \ge 3$. Then merge these contaminated subsets into a new set. Since the new set contains more than three defective items, apply procedure $Z^u$ to it. In this case, the total number of items incurred by procedure $Z^u$ is no more than $n_{\rm \Rmnum{2}}/2$.  The total number of tests in Part~\Rmnum{2} is the sum of four tests for $S^1$, four tests for $S^2$, no more than two tests for $R_2$, and the number of tests incurred by procedure $Z^u$. Thus, based on $3\le d_{\rm \Rmnum{2}}\le n_{\rm \Rmnum{2}}/2$, 
  \begin{eqnarray*}
    t_{\rm \Rmnum{2}} \le 1.431d_{\rm \Rmnum{2}}\big(\log \frac{n_{\rm \Rmnum{2}}/2}{d_{\rm \Rmnum{2}}} + 1.1242\big) + 23 + 4 + 4 + 2 < 1.431d_{\rm \Rmnum{2}}\big(\log \frac{n_{\rm \Rmnum{2}}}{d_{\rm \Rmnum{2}}} + 1.1242\big) + 29,
  \end{eqnarray*}
which holds by Lemma~\ref{lem:UpperBound1-Zu}, and
  \begin{eqnarray*}
  t_{\rm\Rmnum{2}} &\le& 1.4\cdot \frac{n_{\rm \Rmnum{2}}}{2} + 4 + 4 + 2 ~\le~ 1.4n_{\rm \Rmnum{2}} + 6,
  \end{eqnarray*}
which holds by Lemma~\ref{lem:UpperBound2-Zu}.

Based on the above two cases on $\alpha^2$, we have the following results for $\alpha^1 = 2$:
\begin{eqnarray}
  t_{\rm \Rmnum{2}} &\le& 1.431d_{\rm \Rmnum{2}}\big(\log \frac{n_{\rm \Rmnum{2}}}{d_{\rm \Rmnum{2}}} + 1.1242\big) + 29,\label{Part2-2} \\
  {\rm and~}t_{\rm\Rmnum{2}} &<& 1.4n_{\rm \Rmnum{2}} + 13.\label{Part2-2-1}
\end{eqnarray}

If $3\le \alpha^1 \le 4$, then there exist more than three contaminated subsets in $S_1$. Merge them into a new set and apply procedure $Z^u$ to it. Thus, the total number of tests in Part~\Rmnum{2} is four tests for $S_1$ plus the number of tests incurred by procedure $Z^u$. By Lemma~\ref{lem:UpperBound1-Zu}, we have
\begin{eqnarray}
  t_{\rm \Rmnum{2}} &\le& 1.431d_{\rm \Rmnum{2}}\big(\log \frac{n_{\rm \Rmnum{2}}}{d_{\rm \Rmnum{2}}} + 1.1242\big) + 27,\label{Part2-3}
\end{eqnarray}
and based on Lemma~\ref{lem:UpperBound2-Zu},
\begin{eqnarray}
  t_{\rm \Rmnum{2}} \le 1.4n_{\rm \Rmnum{2}} + 4.\label{Part2-3-1}
\end{eqnarray}

From the above analysis of Part~\Rmnum{1} and Part~\Rmnum{2}, we can derive two upper bounds on the number of tests performed by Algorithm $Z^{c}$, which will be used to calculate the final result. Based on formula~(\ref{totaltestZD}), inequalities~(\ref{Part1}),~(\ref{Part2-0}),~(\ref{Part2-1}),(\ref{Part2-2}), and (\ref{Part2-3}), we have the first upper bound such as 
\begin{eqnarray}
  M_{Z^{c}}(d|n) &=& t_{Z^d}\nonumber\\
  &=& t_{\rm\Rmnum{1}} + t_{\rm\Rmnum{2}}\nonumber\\
  &\le& \Big(1.431d_{\rm\Rmnum{1}}(\log \frac{n_{\rm\Rmnum{1}}}{d_{\rm\Rmnum{1}}} + 1.1242) + 3 \Big) + \Big( 1.431d_{\rm \Rmnum{2}}(\log \frac{n_{\rm \Rmnum{2}}}{d_{\rm \Rmnum{2}}} + 1.1242) + 29\Big)\nonumber\\
  &=& 1.431d\big(\log \frac{n}{d} + 1.1242\big) + 32.\nonumber
\end{eqnarray}
Likewise, based on inequalities~(\ref{Part1-1}),~(\ref{Part2-0-1}),~(\ref{Part2-1-1}),~(\ref{Part2-2-1}), and~(\ref{Part2-3-1}), we give the second upper bound such as $M_{Z^{c}}(d|n) \le 1.4n + 13$ for $0\le d \le n$. To sum up, the proof is completed. \hfill$\Box$

\

The total number of tests used by Algorithm~$Z^{c}$ in the worst case for the CGT-U has been given in Lemma~\ref{lem:UpperBound-Zd}. Before analyzing the competitive ratio, we need the following results: Theorem~\ref{theo:I(n-1)} presented by \cite{1982Minimizing}, and Lemma~\ref{lem:1.096} extended from \cite{allemann2003improved}.

\begin{theorem}\label{theo:I(n-1)}
\citep{1982Minimizing}
If $(8/21)n\le d < n$, then $M(d,n)=n-1$.
\end{theorem}

\begin{lemma}\label{lem:1.096}
Let $r=n/d$. For $0< d \le (1/2)n$,
\begin{eqnarray}
M(d,n)\ge d\Big(\log r+(r-1)\log(\frac{r}{r-1})\Big)-0.5\log d-1.5.\label{lemma-8}
\end{eqnarray}
For $r > 21/8$, we have $(r-1)\log(\frac{r}{r-1}) > 1.1243$.
\end{lemma}

{\noindent\bf Proof.}\,~Inequality~(\ref{lemma-8}) is a bound for $\log\binom{n}{d}$ by Stirling's formula and we refer the readers to \cite{allemann2003improved} for a detailed proof. Since $(r-1)\log(\frac{r}{r-1})$ is a monotonically increasing function, we have $(r-1)\log(\frac{r}{r-1}) > 1.1243$ for $r> 21/8$. \hfill$\Box$

\subsection{Proof of Theorem~\ref{theo:paper}}

We are now ready to establish Theorem~\ref{theo:paper} on the competitive ratio of Algorithm~$Z^c$, drawing on Theorem~\ref{theo:I(n-1)} and Lemmas~\ref{lem:1.096}, \ref{lem:UpperBound-Zd}. The proof proceeds by distinguishing the following three cases according to the value of $d$.
  
If $d=0$, according to Step (3) and Step (6), we have $M_{Z^{c}}(d|n)\le 7$.
  
If $(8/21)n\le d < n$, then by Lemma~\ref{lem:UpperBound-Zd} and Theorem~\ref{theo:I(n-1)}, we have
\[
  M_{Z^{c}}(d|n) \le 1.4n+13 < 1.431(n-1)+15 = 1.431M(d,n)+15.
\]
  
If $1\le d < (8/21)n$. Lemma~\ref{lem:UpperBound-Zd} and Lemma~\ref{lem:1.096} yield the following inequalities:
\[
M(d,n) > d(\log \frac{n}{d}+1.1243)-0.5\log d-1.5,
\]
and
\[
M_{Z^{c}}(d|n) \le 1.431d(\log \frac{n}{d}+1.1242)+32.
\]
It follows that
\[
  M_{Z^{c}}(d|n)\le 1.431M(d,n)+32+1.431\delta (d),
\]
where $\delta (d)=0.5\log d +1.5-0.0001d$. By elementary calculus, $\delta(d)$ is an unimodal function for $d\ge 1$ with a unique maximum at $d=500/ \ln 2 \approx 721.35$. For integer $d$, thus,
\begin{eqnarray*}
  M_{Z^{c}}(d|n) &\le& 1.431M(d,n)+32+ 1.431\cdot \max\{\delta(721),\delta(722)\}\\
  &<& 1.431M(d,n)+39
\end{eqnarray*}

Therefore, Theorem~\ref{theo:paper} can be established, implying that Algorithm~$Z^c$ achieves a competitive ratio of $1.431$. In addition, following the analysis of Theorem~\ref{theo:paper}, it is easy to check that the up-zig-zag strategy can achieve a competitive ratio of $1.431$ in the case of $d=0$ and $3 \le d < n$.

\section{Proof of Lemma~\ref{lem:UpperBound1-Zu} and Lemma~\ref{lem:UpperBound2-Zu}}\label{sec:upzigzag}

In this section, we will prove Lemma~\ref{lem:UpperBound1-Zu} and Lemma~\ref{lem:UpperBound2-Zu}. For these, we need to analyze procedure $Z^u$ and bound the performance on the number of tests. To maintain brevity, all steps mentioned in this section are part of Procedure~\ref{alg:Zu}, unless explicitly stated otherwise.

\subsection{Notations, Remark, and Definitions}\label{subsec:NRD}

We begin with a remark, several notations, and definitions, some of which are already used in \cite{wu2022improved}. Here, we restate them for the sake of completeness.

First, we redefine the status of a set in Remark~\ref{rem:redefine}. 
\begin{remark}\rm
  \label{rem:redefine}
  In general, we say a set is \emph{contaminated} while there exists at least one defective item in it. However, for a set $S'$ of size $2$, we say that it is ‘\emph{contaminated}’ only if $S'=\{1,1\}$, else is ‘\emph{pure}’ during the process of procedure $Z^u$.
\end{remark}

Define $n(Z^u)$ and $d_{Z^u}$ as the total number of items and the number of defective items identified by procedure $Z^u$, respectively. As shown in Procedure~\ref{alg:Zu}, procedure $Z^u$ contains two types of tests: (\rmnum{1}) the first type consists of tests performed in Steps (5) and (13); (\rmnum{2}) the second type consists of tests performed in Steps (20), (24) and (31). We use $\mathscr{T} = \{ T_1, T_2,\cdots, T_q \}$ to denote a set of all the tests of the first type. In order to delineate the process of procedure $Z^u$, an \emph{additional test} $T'$ in $\mathscr{T}$ and the concept of \emph{phase} are introduced.

\begin{definition}\rm
  \label{def:additional}
  {\sc (Additional Test)} We call $T'\in \mathscr{T}$ an \emph{additional test} only if $T'$ is a test performed in Step (5) of Procedure~\ref{alg:Zu}. Thus, any additional test before $T_q$ is contaminated. 
\end{definition}

\begin{definition}\rm
  \label{def:phase}
  {\sc (Phase)}~We refer to such a testing sequence $\mathscr{P} = \{T_i, T_{i+1}, \cdots, T_{i+j}\}\subseteq \mathscr{T}$ for $j\ge 0$ as a ‘\emph{phase}’ if it satisfies one of the following four conditions:
  \begin{itemize}
    \item when $j = 0$, the test result of $T_{i+j}$($ = T_i$) is contaminated and the test result of $T_{i-1}$ is contaminated if it exists;
    \item when $1\le j\le 5$, the test result of $T_{i+j}$ is contaminated and the others' in $\mathscr{P}$ are pure;
    \item when $j\ge 6$, the test results of $T_{i+6}$ and $T_{i+j}$ are contaminated, and the others' in $\mathscr{P}$ are pure. Here, $T_{i+6} = T_{i+j}$ if $j=6$; otherwise, $T_{i+6}$ is an additional test. 
    \item when $i+j = q$, the test results of ${T_i, T_{i+1}, \cdots, T_{i+j}}$ are all pure and the test result of $T_{i-1}$ is contaminated if it exists.
  \end{itemize}
\end{definition}

Namely, a \emph{phase} is a process that is from the end of identifying the preceding defective item (or the first test at the beginning of procedure $Z^u$, i.e., $T_1$) to the end of identifying at least one defective item from the current contaminated set (or identifying that all remaining items are good). Thus, a path of procedure $Z^u$ consists of no more than $d_{Z^u}+1$ phases. Let $\mathscr{P}_i$ be the set of tests in the $i^{th}$ phase.

The \emph{rank} of a test in $\mathscr{T}$ and \emph{incurred tests} are defined.

\begin{definition}\rm
  \label{def:rank}
  {\sc (Rank)}~During the execution of procedure $Z^u$, excluding all additional tests, for any other test $T$ in $\mathscr{T}$ (i.e., performed in Step~(13)), if the current value of $k$ is $v$, we say that $T$ is of \emph{rank} $v$, defined as $T^v$.
\end{definition}

The values of ranks in a phase are successive, or it is exactly equal to zero when there is only one test in the phase. For any test $T\in \mathscr{T}$, let $r(T)$ be the rank of $T$.

\begin{definition}
  \label{def:incurredtests}\rm
  {\sc (Incurred Tests)}~For a test $T\in \mathscr{T}$, let $S$ be the tested set of $T$. Define $I(T)$ as the set of tests \emph{incurred} by $T$ as follows:
  \begin{itemize}
    \item if the test result of $T$ is pure, then $I(T)$ contains only one test, i.e., $I(T) = \{T\}$;
    \item if $T$ is an additional test, then $I(T)$ contains only one test, i.e., $I(T) = \{T\}$;
    \item if the test result of $T$ is contaminated and not an additional test, then $I(T)$ contains all tests performed on $S$.
  \end{itemize} 
\end{definition}

By the above definitions, the set of all tests performed by procedure $Z^u$ can be written as $\bigcup_{T\in \mathscr{T}}I(T)$, and the total number of tests used by procedure $Z^u$ is
\begin{equation}\label{Zutests}
  t_{Z^u} = \sum_{T\in \mathscr{T}}|I(T)| = \sum^q_{j = 1}|I(T_j)|.
\end{equation}

Next, the \emph{identified items} of a test $T$ and a \emph{zig-zag tuple} consisting of two or three tests in $\mathscr{T}$ should be introduced.

\begin{definition}\rm
  \label{def:identified items}
  {\sc (Identified Items)}~For a test $T\in \mathscr{T}$, define $D(T)$ to be the set of items \emph{identified} by tests in $I(T)$, and define $n(T) = |D(T)|$.
\end{definition}

According to procedure $Z^u$, for any nonnegative integer $v$, a test of rank $v$ is applied to a subset of size at most $a_v$.

\begin{definition}\rm
  \label{def:zigzagtuple}
  {\sc (Zig-Zag Tuple)}~We say that a \emph{zig-zag tuple} is form of three tests $T_i, T_j$, and $T'_k$ in $\mathscr{T}$ such that $P = (T_i, T_j, T'_k)$. For the first two tests in $P$, $T_i$ has rank $v-1$ with result pure, and $T_j$ has rank $v$ with result contaminated, where $v\ge 1$, moreover $T_i$ is required to be applied to a subset of size exactly $a_{v-1}$. The final test $T'_k$ is an additional test that is in the same phase $\mathscr{P}'$ as $T_j$. Note that only if $\mathscr{P}'$ contains no less than six tests with pure results and $|D(T_j)| = a_v$, an additional test would be occurred, thus, $|I(T'_k)|=1$; otherwise, $I(T_k') = \emptyset$, thus $|I(T'_k)|=0$, and we write this kind of tuple in a simple form $P = (T_i, T_j, \emptyset)$. For brevity, let $v$ be the rank of tuple $P$.
\end{definition}

In the definition of a zig-zag tuple, $P = (T_i, T_j, T'_k)$, $T_i$ and $T_j$ are not required to be in the same phase (consecutive) in $\mathscr{T}$, and there is no requirement on the order of $T_i$ and $T_j$, that is, both $i < j$ and $i > j$ are allowed. In addition, we have $k \le j - 1$. Thereafter, we will simply refer to a ‘zig-zag tuple’ as a ‘tuple’ for short if there is no confusion.

\subsection{Partition {$\mathscr{T}$} Into Four Classes}\label{subsec:Partition}

\begin{figure}[!bph]
  \centering
   \includegraphics[width=0.85\textwidth]{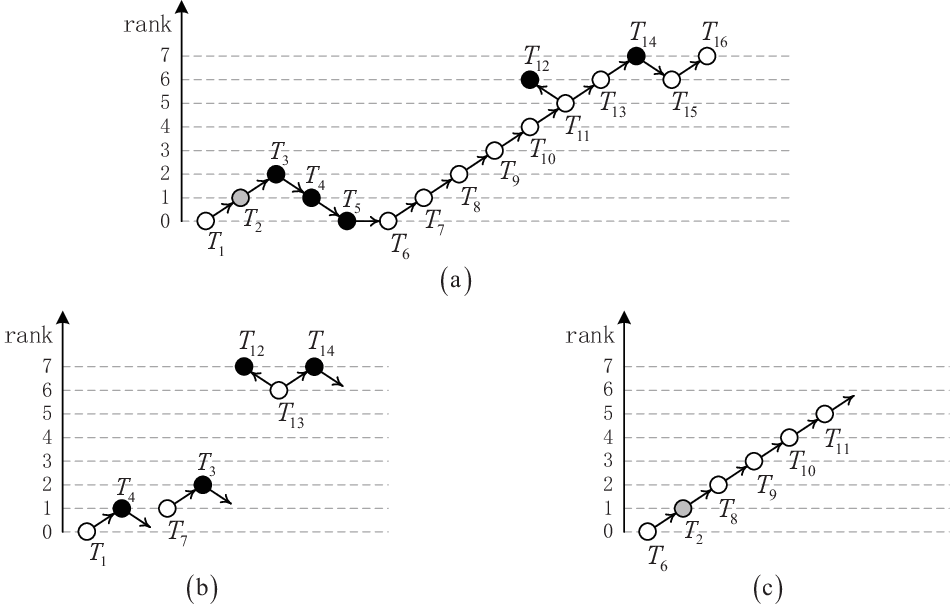}
   \caption{.~{\bf An Example to Illustrate the Partition of $\mathscr{T}$ into Four Classes.} \\ 
   \emph{Note.}~{\bf(a)}~The circles represent tests performed in Steps (5) and (13) (i.e., tests in $\mathscr{T}$), in order that they are performed by Procedure~\ref{alg:Zu}. In this example, we have $\mathscr{T} = \{T_1, T_2,\cdots, T_{16}\}$. This test process contains $5$ phases, such as $\mathscr{P}_1 = \{T_1, T_2, T_3\}$, $\mathscr{P}_2 = \{T_4\}$, $\mathscr{P}_3 = \{T_5\}$, $\mathscr{P}_4 = \{T_6, T_7, T_8, T_9, T_{10}, T_{11}, T_{12}, T_{13}, T_{14}\}$ and $\mathscr{P}_5 = \{T_{15}, T_{16}\}$. An arrow pointing to the bottom right from a circle (test) $T$ indicates that the test result of
   $T$ is contaminated, and the rank of the test (performed in Step (13)) immediately following $T$ is decreased by one. An arrow pointing
   to top right from a test $T$ indicates that the test result of $T$ is pure, and if $T$ is the sixth test in this phase, the test (performed in Step (5)) immediately following $T$ is for the whole remaining items; otherwise, the rank of the test (performed in Step (13)) immediately following $T$ is increased by one. An arrow pointing to the top left from a test $T$ indicates that the test result of $T$ is contaminated and $T$ is an \emph{additional test}. Assume the preceding test of an additional test is $T'$, then the rank of the test (performed in Step (13)) immediately following an additional test is increased by one based on $r(T')$, i.e., $r(T') + 1$. Notice that a test $T$ of rank zero may have a contaminated test result. In this case, a horizontal arrow is drawn from $T$ to the right, indicating that the next test in Step (13) also has rank zero. For this example, $C_1 = \{T_{15}, T_{16}\}$ with $\Delta_{C_1}=1$, and $C_2=\{T_5\}$. {\bf(b)}~Tests in $C_3$. In this example, $C_3$ is formed by three zig-zag tuples, and one possibility for the tuples could be $(T_1, T_4,\emptyset)$, $(T_7, T_3,\emptyset)$ and $(T_{13}, T_{14}, T_{12})$. Here $T_{12}$ is an additional test with $D(T_{12}) = D(T_{13})\cup D(T_{14})\cup D(T_{15})\cup D(T_{16})$. {\bf(c)}~Tests in $C_4 = \mathscr{T}\setminus (C_1\cup C_2\cup C_3) = \{T_6, T_2, T_8, T_9, T_{10}, T_{11}\}$. Notice that $T_2$ (gray circles) is a \emph{pure} test for a subset \{0,1\} and $T_7$ is a \emph{pure} test for a subset $\{0,0\}$ based on Remark~\ref{rem:redefine}. So $C_3$ is not unique, for example, we can also have $C_3 = \{(T_1, T_4,\emptyset), (T_2, T_3,\emptyset),(T_{13}, T_{14}, T_{12})\}$.}
   \label{fig:example}
 \end{figure}

We partition $\mathscr{T} = \{T_1, T_2, \cdots, T_q\}$, the set of all tests performed by procedure $Z^u$, into four classes $C_1$, $C_2$, $C_3$, and $C_4$. These four classes are defined as follows (an example is shown in Figure~\ref{fig:example}):
\begin{itemize}
  \item~$C_1$ is formed by all the tests in the final phase, i.e., several consecutive tests defined as $T_{q-\Delta_{C_1}}, T_{q-\Delta_{C_1}+1}, \cdots, T_q$ with results pure, where $\Delta_{C_1}$ is a nonnegative integer. Once there exist more than five consecutive pure subsets in the final phase, the remaining items after the fifth test would be tested as a single entity, so $0\le \Delta_{C_1} \le 6$;
  \item~$C_2$ is formed by all tests having a contaminated result and rank zero;
  \item~$C_3$ is formed by three categories of tests. The first category refers to all additional tests. The second category refers to every test $T$ having a result contaminated in $\mathscr{T}\setminus C_2$. For each $T$ in the second category, we find a test having rank $r(T)-1$ and result pure from $\mathscr{T}$ corresponding to it, and then group this test into the third category;
  \item~$C_4$ is formed by all the tests in $\mathscr{T}\setminus (C_1\cup C_2\cup C_3)$.
\end{itemize}

\

Subsequently, we have the following two claims.

\begin{claim}
  \label{ob:ComPath}~Suppose that $A$ is a set of tests formed by all additional tests in $\mathscr{T}$, then tests in $C_3$ can be composed into $(|C_3|- |A|)/2$ zig-zag tuples, such that each test in $C_3$ appears in exactly one tuple.
\end{claim}

\noindent {\bf Proof.}\,~From the definition of class $C_1$, class $C_2$, and additional tests, after all the tests belonging to them are removed from $\mathscr{T}$, the ranks of any two consecutive tests in the remaining subset $\mathscr{T}'$ differ by one (note that by ‘consecutive’ we refer to the sequence after removal of $C_1\cup C_2$ and all the additional tests, that is, two tests could be consecutive on rank in the remaining test set $\mathscr{T}' \in \mathscr{T}$). Indeed, the final test on $\mathscr{T}'$ has a contaminated test result. We assume that $r_m$ is the largest rank of a test with a result contaminated in $\mathscr{T}'$.

Consider a vertical line segment $L$ (which is part of the $y$-axis) with a lower endpoint $y=0$ and an upper endpoint $y=r_m$. We map each test $T\in \mathscr{T'}$ to a \emph{move} between two consecutive integer points in $\{y=0,y=1,\cdots,y=r_m\}$ of $L$ as follows. The mapping depends on the test result of $T$. If the test result of $T_x$ is pure, we have $0\le r(T_x)\le r_m - 1$, because if $r(T_x)\ge r_m$ then the test result of some test next $T$ must be contaminated, contradicting the assumption of $r_m$. We map $T_x$ to a move on $L$ from point $y=r(T_x)$ to point $y=r(T_x)+1$, and use $(r(T_x)\rightarrow)$ to denote the move. If the test result of $T_y$ is contaminated, then since $T_y\notin C_2$, we have $1\le r(T_y)\le r_m$, and we map $T_y$ to a move on $L$ from point $y=r(T_y)$ to point $y=r(T_y)-1$ and use $(r(T)\Rightarrow)$ to denote the move. By the above mapping, the tests in $\mathscr{T}'$ map to a walk $W^*$, which is a set consisting of a sequence of moves, and each move may either go upwards or go downwards on $L$. Since the ranks of the tests in $\mathscr{T}'$ are consecutive, there must exist a one-to-one correspondence between a move $(r(T_x)\rightarrow)$ and a move $(r(T_y)\Rightarrow)$ where $r(T_x)=r(T_y)-1$, if not then contradicts the step of the rank going downwards or upwards by 1 once a time. Thus, we can find as many moves $(r(T)-1\rightarrow)$ as moves $(r(T)\Rightarrow)$ in $W^*$.

In addition, by Step (4) and Step (5), it will occur an additional test $T^{'}_{z}\in A \subseteq \mathscr{T}$, where $T^{'}_{z}$ is a succeeding test of $T_x$ with result contaminated, only if $T_x$ is the sixth test with result pure in the current phase. Given that each $T^{'}_{z}$ is of result contaminated, it follows that there must exist a test (such as $T_y$) in this phase. Consequently, the number of additional tests similar to $T^{'}_{z}$ is equal to or less than the number of tests similar to $T_y$. 

From the above analysis and by the definition of a zig-zag tuple, for any $1\le r(T)\le r_m$, the two tests mapping to the moves $(r(T)-1\rightarrow)$ and $(r(T)\Rightarrow)$ form a zig-zag tuple with a possible additional test $T'$ (note that $T'$ exists only if it occurs more than five tests with results pure in the phase containing $T$), implying that the number of tuples in $C_3$ equals to $(|C_3|-|A|)/2$, which establishes the claim (see Figure~\ref{fig:example}(b)). \hfill $\Box$

\begin{claim}
  \label{ob:rank}
  Suppose that a test $T'$ whose result is contaminated has the largest rank in $C_3$, and that $T''$ has the largest rank in $C_4$. Then $C_4$ consists of a series of disjoint tests with consecutive ranks $(0, 1, 2, \cdots, r(T''))$, all of which have pure results, where $r(T'')+2 \le r(T')$.
\end{claim}

\noindent {\bf Proof.}\,~We first prove that $C_4$ is formed by a series of disjoint tests of consecutive ranks such that $(0,1,2,\cdots,r(T''))$ if it is not empty.

Since $C_4 = \mathscr{T} \setminus (C_1 \cup C_2 \cup C_3)$, it is evident that every test in $C_4$ is indeed pure in terms of its test results. Consider the case of $C_4 = \{T_u, T_v, \cdots, T''\}$ (note that we rearrange the tests in $C_4$ based on their ranks in ascending order), for the sake of contradiction, we assume that $r(T_u)+1 < r(T_v)$. The ranks of the tests in the subset $\mathscr{T'}$(consisting of the remaining tests after removing all additional tests from $\mathscr{T}\setminus (C_1\cup C_2)$) are consecutive, so there exists a test subset $\mathscr{T}_1$ between $T_u$ and $T_v$, such that it contains $r(T_v)-r(T_u)-1$ pure tests of ranks $(r(T_u)+1,r(T_u)+2, \cdots, r(T_v)-1)$. Since $\mathscr{T}_1$ does not belong to $C_4$ and is exactly not in $C_1\cup C_2$, we have $\mathscr{T}_1\in C_3$. By the definition of $C_3$, it has another test subset, denoted by $\mathscr{T}_2$, comprising a consecutive series of contaminated tests with ranks $(r(T_u)+2,r(T_u)+3, \cdots, r(T_v))$, where each test in $\mathscr{T}_2$ is paired individually with every test in $\mathscr{T}_1$. Let $T_{z}$ be the test with contaminated result and rank $r(T_v)$ in $\mathscr{T}_1$. Based on the continuity of ranks, there must exist a test subset denoted $\mathscr{T}_3$, consisting of a series of tests with pure results and ranks $(r(T_u), r(T_u) + 1, \cdots, r(T_v))$, positioned between $T_z$ and $T_v$. However, in this way, $\mathscr{T}_3 \subseteq C_4$ which contradicts that there is no test between $T_u$ and $T_v$. As a result, we can conclude that $r(T_u)+1 \ge r(T_v)$. 

Next, we establish the proof by contradiction to demonstrate that $T_v$ is the sole test in $C_4$ with a rank of $r(T_v)$. Assume that there exists another test $T_w$ with rank $r(T_w) = r(T_v)$ after $T_v$, i.e., $w > v$. Based on this assumption, $T_{w'}$ (the preceding test of $T_w$, $T_{w'}\in \mathscr{T}'$) has result contaminated and $r(T_{w'}) = r(T_w)+1 = r(T_v)+1$. According to Claim~\ref{ob:ComPath}, $\{T_v, T_{w'}\} \in C_3$, which contradicts $T_v \in C_4$. So $T_{w}$ does not exist and $T_v$ is the only test with rank $r(T_v)$ in $C_4$, which means $r(T_u) < r(T_v)$. Combined with $r(T_u)+1 \ge r(T_v)$, we can get $r(T_u)+1=r(T_v)$. In addition, procedure $Z^u$ begins with a test of rank $0$, and the rank of any test $T\in \mathscr{T}$ is nonnegative, so $r(T_u) = 0$. As such, $C_4$ consists of a sequence of disjoint tests having consecutive ranks $(0,1,2,\cdots,r(T''))$ and results pure. Additionally, for the cases where $C_4 = \emptyset$ or $C_4 = \{T''\}$, it is evident that the claim hold.

Since $T''\notin C_1$, there exists at least a test $T_x$ of test result pure and a test $T_y$ of test result contaminated after $T''$. Thus, $r(T'')+2 \le r(T_x)+1 \le r(T_y)$. Moreover, $T'$ is a test with the largest rank in $C_3$, so $r(T') \ge r(T_y)$. Hence, we have $r(T'')+2 \le r(T')$. Therefore, the claim is proven. \hfill $\Box$

\subsection{Upper Bounds on Classes {$C_1$}, {$C_2$}, {$C_3$}, and {$C_4$}}\label{subsec:
Estimating}

We have partitioned $\mathscr{T}$ into four disjoint subsets $C_1$, $C_2$, $C_3$, and $C_4$. Thus, by equation~(\ref{Zutests}), the total number of tests used by procedure $Z^u$ is
\begin{equation}\label{Zualltests}
  t_{Z^u}=\sum_{T\in \mathscr{T}} |I(T)|=\sum_{T\in C_1} |I(T)|+\sum_{T\in C_2} |I(T)|+\sum_{T\in C_3} |I(T)|+\sum_{T\in C_4}|I(T)|.
\end{equation}

For formula~(\ref{Zualltests}), we use a technique so-called divide-and-conquer to analyze: first estimate the first term, the second term, and the third term separately; then estimate the sum of $\sum_{T\in C_2\cup C_3\cup C_4} |I(T)|$; finally, incorporate the first term to reach an upper bound on $t_{Z^u}$.

\begin{lemma}\label{C_1}
  For $C_1$, we have $\sum_{T\in C_1}|I(T)| \le 7$.
\end{lemma}

\noindent{\bf Proof.}\, By Step $(4)$ and the definition of $C_1$, it is easy to obtain that the number of all tests in $C_1$ does not exceed $7$. \hfill $\Box$

\begin{lemma}\label{C_2}
  For each test $T\in C_2$, we have $|I(T)| < 1.431\big(\log n(T)+1.1242\big)$.
\end{lemma}

\noindent {\bf Proof.}\, By the definition of $C_2$, each test $T\in C_2$ is of rank 0 and test result contaminated. By Step $(20)$ in procedure 4-Split, $T$ is applied to a subset of size one, and it identifies the only item as defective. Thus, $n(T)=1$, $I(T)={T}$ and $|I(T)|=1$. Therefore, we have $|I(T)| = 1 < 1.431\big(\log n(T)+1.1242\big)$ for each $T\in C_2$. \hfill $\Box$

\

According to Claim~\ref{ob:ComPath}, all tests in $C_3$ can be divided into $|C_3-\sum T'_k|/2$ tuples, such that we can calculate the total test number based on each tuple. Let $C_3P$ denote the set of all tuples obtained from $C_3$. So for any tuple $P=(T_a,T_b,T'_c)$ in $C_3P$, $I(P) = I(T_a)\cup I(T_b)\cup I(T'_c)$ based on Definition~\ref{def:identified items}. Let $D(P)$ denote the set of all items identified by tests in $I(P)$, thus, $D(P)=D(T_a)\cup D(T_b)$ and $n(P)=|D(P)|$. Thus,
$$
|I(P)|=|I(T_a)|+|I(T_b)|+|I(T'_c)|,
$$
$$
n(P)=|D(P)|=|D(T_a)|+|D(T_b)|=n(T_a)+n(T_b).
$$
Subsequently, the second term in equation~(\ref{Zualltests}) can be written as
$$
\sum_{T\in C_3} |I(T)|=\sum_{T\in C_3P} |I(P)|.
$$

Before estimating the number of tests in $C_3P$, we first refer to some special tuples based on Remark~\ref{rem:redefine}. Remark~\ref{rem:specialtuple} presents descriptions of various noteworthy tests, while Claim~\ref{ob:specialtuple} enumerates five specific tuples, highlighting their significance within the context of the study.

\begin{remark}\rm
  \label{rem:specialtuple}
  Given a test, $T^v$, performed on a subset of size $a_v$. Let $^dT^v$ be one type of test on $T^v$, where $d$ is the number of defective items identified by the test.
  \begin{itemize}
    \item For $T^1$, there exist three types of tests depending on $d$:~{\rm(\rmnum{1})}~$^0T^1$: $D(^0T^1)=\{0,0\}$, $|I(^0T^1)| = 1$; {\rm(\rmnum{2})}~$^1T^1$: $D(^1T^1)=\{0,1\}$, $|I(^1T^1)| = 3$; {\rm(\rmnum{3})}~$^2T^1$: $D(^2T^1)=\{1,1\}$, $|I(^2T^1)| = 3$.
    \item For $T^2$ (it is contaminated and $\Delta_{\rm TS}$ in the current phase equals $1$), there exist three types of tests depending on $d$:~{\rm(\rmnum{1})}~$^1T^2$: $D(^1T^2)=\{0,0,1\}$, $|I(^1T^2)| = 4$;~{\rm(\rmnum{2})}~$^2T^2$: $D(^2T^2)=\{0,1,1\}$, $|I(^2T^2)| = 3$;~{\rm(\rmnum{3
    })}~$^3T^2$: $D(^3T^2)=\{1,1,1\}$, $|I(^3T^2)| = 4$.
  \end{itemize} 
\end{remark}

\begin{figure}[ht]
  \centering
   \includegraphics[width=0.6\textwidth]{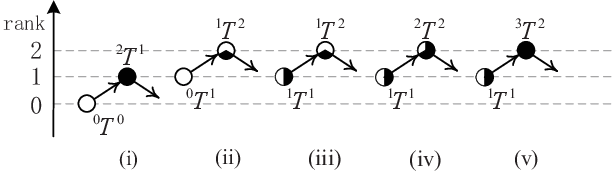}
   \caption{.~{\bf An Example to Illustrate the Five Types of Tuples in Claim~\ref{ob:specialtuple}} \\
   \emph{Note.}~The serial numbers \{(\rmnum{1}),$\cdots$,(\rmnum{5})\} in the figure correspond to the serial numbers in Claim~\ref{ob:specialtuple}. In order to distinguish the different tests on $T^1$ and $T^2$, we use the percentage of black areas in a hollow circle to mark the number of defective items detected by a test. For example, $D(^2T^2) = \{0,1,1\}$ and so $^2T^2$ is a cycle covered by two-thirds of the black area.}
   \label{fig:special tuples}
\end{figure}

\begin{claim}
  \label{ob:specialtuple}
  For each zig zag tuple $P\in C_3P$, there exist the following five possible types of tuples on rank 1 and 2, where $d_p$ is the number of defective items detected by $P$. Refer to the illustration in Figure~\ref{fig:special tuples}.
  \begin{itemize}\rm
    \item[(\rmnum{1})]~$P=(^0T^0,^2T^1,\emptyset)$: $|I(P)|=4$, $n(P)=3$, $d_P=2$.
    \item[(\rmnum{2})]~$P=(^0T^1,^1T^2,\emptyset)$: $|I(P)|=3({\rm or}~ 4)$, $n(P)=3(4~{\rm or}~5)$, $d_P=1$.
    \item[(\rmnum{3})]~$P=(^1T^1,^1T^2,\emptyset)$: $|I(P)|=7$, $n(P)=5$, $d_P=2$.
    \item[(\rmnum{4})]~$P=(^1T^1,^2T^2,\emptyset)$: $|I(P)|=7$, $n(P)=5$, $d_P=3$.
    \item[(\rmnum{5})]~$P=(^1T^1,^3T^2,\emptyset)$: $|I(P)|=7$, $n(P)=5$, $d_P=4$.
  \end{itemize}
\end{claim}

For the sake of description, we hereby make a trivial claim about a tuple $P_x$.
\begin{claim}
  \label{ob:A}
  For a zig zag tuple $P_x = \{T_x^{i_x-1}, T_x^{i_x}, T'_x\}$ of rank $i_x$ where $i_x\ge 3$, procedure 4-Split is applied to $D(T_x^{i_x})$. Partition $T_x^{i_x}$ into four subsets $Y_x, Z_x, U_x$, and $V_x$, and so there are four types of test results as follows:
  \begin{itemize}\label{P_x}
    \item {\rm\bf T1.}~$Y_x$ is contaminated. $|I(P_x)| \le |I(T'_x)|+i_x+1$, $n(P_x) \ge 3\cdot 2^{i_x-3} + 1$;
    \item {\rm\bf T2.}~$Y_x$ is pure and $Z_x$ is contaminated. $|I(P_x)| \le |I(T'_x)|+ i_x+2$, $n(P_x) \ge 5\cdot 2^{i_x-3} + 1$;
    \item {\rm\bf T3.}~$Y_x$ and $Z_x$ are pure, $U_x$ is contaminated. $|I(P_x)| \le |I(T'_x)|+ i_x+2$, $n(P_x) \ge 7\cdot 2^{i_x-3} + 1$;
    \item {\rm\bf T4.}~$Y_x$, $Z_x$, and $U_x$ are pure, $V_x$ is contaminated. $|I(P_x)| \le |I(T'_x)|+ i_x+2$, $n(P_x) \ge 8\cdot 2^{i_x-3} + 1$.
  \end{itemize}
\end{claim}

\noindent{\bf Proof Sketch.}\,~The proof is trivial and so we put it into Appendix~\ref{sec:appendix1}. \hfill $\Box$

\

Next, we give a bound on each tuple as follows.

\begin{lemma}
  \label{C_3}
  For each zig-zag tuple $P= (T_a, T_b, T'_c)$, we have 
  $$
  |I(P)|\le 1.431d_P\big(\log \frac{n(P)}{d_P}+1.1242\big),
  $$
  where $d_P$ is the number of defective items identified by $P$.
\end{lemma}

\noindent{\bf Proof.}\,~Let $i$ be the rank of $P= (T_a, T_b, T'_c)$. By the definition of $C_3$, we have $i\ge 1$, and $T_a$ has result pure and is applied to a subset $S_a$ of size exactly $3\cdot 2^{i-3}$; $T_b$ has test result contaminated, and is applied to (together with Steps (18),(28) or (32) incurred by $T_b$) a subset $S_b$ of size at most $3 \cdot 2^{i-2}$; $T'_c$ has test result contaminated while $T'_c \neq \emptyset$. Thus, a certain number of defective items and an unspecified number of good items from $S_b$ are identified, and all items in $S_a$ are identified as good items. It follows that the total number of identified items satisfies the following cases.

If $i = 1$, then $|S_a| = 1$ and $|S_b| \le 2$. Procedure 4-Split performs on items in $S_b$. By Claim~\ref{ob:specialtuple}, $P=(^0T^0,^2T^1,\emptyset)$. The total number of tests $|I(P)|$ incurred by $S_a$ and $S_b$ equals $4$, and the total number of identified items satisfies $n(P) = 3$ and $d_p = 2$, thus,
$$
|I(P)|=4 < 1.431\times 2 \times \big(\log(3/2)+1.1242\big) = 1.431d_P\big(\log \frac{n(P)}{d_P}+1.1242\big).
$$

If $i = 2$, then $|S_a| = 2$ and $|S_b| \le 3$. By Claim~\ref{ob:specialtuple}, it has four kinds of tuples on $P$. We conclude with the following two cases: {Case~1a} and {Case~1b}.

{{\bf Case 1a.}}~$P=(^0T^1,^1T^2,\emptyset)$. By Step (\ref{S'<=3}) of Procedure~\ref{alg:4Split}, procedure 4-Split is applied to $S_b$. According to the items' order in $S_b$, we assume that $S_b:=\{b_1,\cdots,b_t\}$ for $1\le t\le 3$. So all the possible test outcomes are the following:
\begin{itemize} 
  \item Item $b_1$ is contaminated. $I(P) = 3$, $n(P)=3$, $d_P = 1$. Thus,
  $$
  |I(P)|=3 < 1.431\times(\log3+1.1242);
  $$
  \item Item $b_1$ is pure and item $b_2$ is contaminated. $I(P) = 4$, $n(P)=4$, $d_P = 1$. Thus,
  $$
  |I(P)|=4 < 1.431\times(\log4+1.1242);
  $$
  \item Items $b_1$ and $b_2$ are pure, and item $b_3$ is contaminated.
  $I(P) = 4$, $n(P)=5$, $d_P = 1$. Thus,
  $$
  |I(P)|=4 < 1.431\times(\log5+1.1242).
  $$
\end{itemize}

{{\bf Case 1b.}}~$P=(^1T^1,^1T^2\setminus ^2T^2\setminus ^3T^2,\emptyset)$. By Claim~\ref{ob:specialtuple} and Step $(25)$, the total number of tests $|I(P)|$ incurred by $S_a$ and $S_b$ is $7$, $n(P) = 5$, and $d_P$ is equal to $2$,$3$ or $4$.
\begin{itemize}
  \item $d_P = 2$: 
  $
  |I(P)|=7 < 1.431\times 2\times \big(\log(5/2)+1.1242\big);
  $
  \item $d_P = 3$:
  $
  |I(P)|=7 < 1.431\times 3\times \big(\log(5/3)+1.1242\big);
  $
  \item $d_P = 4$:
  $
  |I(P)|=7 < 1.431\times 4\times \big(\log(5/4)+1.1242\big).
  $
\end{itemize}
Based on Case 1a and Case 1b, we have
$
|I(P)|\le 1.431 d_P \big(\log \frac{n(P)}{d_P}+1.1242\big)
$
for $i=2$.

\

For $i\ge 3$, procedure 4-Split is applied to $S_b$, and we simply refer to that $S_b$ will be partitioned into four subsets $Y$, $Z$, $U$, and $V$ in the following. Based on Claim~\ref{ob:A}, we have the following four cases {Case 2a}, {Case 2b}, {Case 2c}, and {Case 2d} corresponding to {\bf T1}, {\bf T2}, {\bf T3}, and {\bf T4}, respectively.

{{\bf Case 2a.}}~$Y$ is contaminated. In this case, we have $|I(P)|\le |I(T'_c)|+i+1$ and $n(P) \ge |S_a|+1= 3\cdot2^{i-3}+1$. If $T_a$ is the sixth test in its phase and $|S_b| = 3\cdot 2^{i-2}$, then $T'_c \neq \emptyset$ and $|I(T'_c)| = 1$; otherwise, $T'_c=\emptyset$ and $|I(T'_c)| = 0$. 

Let $f^1_Y=1.431 d_P \big(\log \frac{3\cdot2^{i-3}+1}{d_P}+1.1242\big)-|I(P)| \le 1.431 d_P \big(\log \frac{n(P)}{d_P}+1.1242\big)-|I(P)|$.
\begin{itemize}
  \item When $|I(T'_c)| = 0$, $|I(P)|=i+1$ for $i \ge 3$,
  \begin{eqnarray*}
  f^1_Y = 1.431\cdot\big(\log(1+3\cdot2^{i-3})+1.1242\big)-i-1.
  \end{eqnarray*}
  Since $f^1_Y$ is a monotonically increasing function for $i\ge 3$, then $f^1_Y\ge f^1_Y(3) >0$.
  \item When $|I(T'_c)| = 1$, $|I(P)|=i+2$ for $i \ge 6$,
  \begin{eqnarray*}
    f^2_Y = f^1_Y - 1 = 1.431\cdot\big(\log(1+3\cdot2^{i-3})+1.1242\big)-i-2.
  \end{eqnarray*}
  Compare to $f^1_Y$, $f^2_Y \ge f^2_Y(6) >0$.
\end{itemize}
Hence, it follows that
$
|I(P)|\le 1.431 d_P \big(\log \frac{n(P)}{d_P}+1.1242\big).
$

{{\bf Case 2b.}}~$Y$ is pure and $Z$ is contaminated. In this case, $|I(P)|\le |I(T'_c)|+i+2$ and $n(P) \ge |S_a|+|Y|+1=5\cdot2^{i-3}+1$. Similar to {Case 1b}, if $T_a$ is the sixth test in the current phase and $|S_b| = 3\cdot 2^{i-2}$, $|I(T'_c)| = 1$; otherwise, $|I(T'_c)| = 0$. Thus, $|I(P)|$ is at most $i+3$.

Let $f^1_Z=1.45 d_P \big(\log \frac{5\cdot2^{i-3}+1}{d_P}+1.1242\big)-|I(P)| \le 1.431 d_P \big(\log \frac{n(P)}{d_P}+1.1242\big)-|I(P)|$.
\begin{itemize}
  \item When $|I(T'_c)| = 0$, $|I(P)|=i+2$ for $i \ge 3$,
  \begin{eqnarray*}
    f^1_Z = 1.431\cdot\big(\log(1+5\cdot2^{i-3})+1.1242\big)-i-2.
  \end{eqnarray*}
  Since $f^1_Z$ is a monotonically increasing function for $i\ge 3$, thus,
  $
  f^1_Z \ge f^1_Z(3) >0.
  $
  \item When $|I(T'_c)| = 1$, $|I(P)|=i+3$ for $i\ge6$,
  \begin{eqnarray*}
    f^2_Z = f^1_Z - 1 \ge 1.431\cdot\big(\log(1+5\cdot2^{i-3})+1.1242\big)-i-3.
  \end{eqnarray*}
  Then similarly, $f^2_Z \ge f^2_Z(6) >0$.
\end{itemize}
Hence, it follows that++
$
|I(P)|\le 1.45 d_P \big(\log \frac{n(P)}{d_P}+1.1242\big).
$

{{\bf Case 2c.}}~$Y$ and $Z$ are pure, $U$ is contaminated. In this case, we have $|I(P)| \le |I(T'_c)|+i+2$ and $n(P) \ge |S_a|+|Y|+|Z|+1=7\cdot2^{i-3}+1$. Compared with {Case 2b}, the number of tests in these two cases is the same, but more items are identified in {Case 2c} than in {Case 2b}. Thus, it follows that
\begin{eqnarray*}
  |I(P)| \le 1.431\big(\log(1+5\cdot2^{i-3})+1.1242\big) <  1.431\big(\log(1+7\cdot2^{i-3})+1.1242\big).
\end{eqnarray*}

{{\bf Case 2d.}}~$Y$, $Z$, and $U$ are pure, $V$ is contaminated. In this case, we have $|I(P)|\le i+2+|I(T'_c)|$ and $n(P) \ge |S_a|+|Y|+|Z|+|U|+1=8\cdot2^{i-3}+1$. Compared with {Case 2c}, the number of tests in these two cases is the same, but more items are identified in {Case 2d} than in {Case 2c}. Thus, it follows that
\begin{eqnarray*}
  |I(P)| < 1.431\big(\log(1+7\cdot2^{i-3})+1.1242\big) < 1.431\big(\log(1+8\cdot2^{i-3})+1.1242\big).
\end{eqnarray*}

To sum up, $I(P)\le 1.431d_P\big(\log \frac{n(P)}{d_P}+1.1242\big)$, as desired. \hfill $\Box$

\

Based on Lemma~\ref{C_2} and Lemma~\ref{C_3}, we have the following upper bound on the sum of the second and third terms in equation~(\ref{Zualltests}).

\begin{lemma}
  \label{C2C3}
  For $C_2\cup C_3$, we have
  \[{\displaystyle \sum_{T\in C_2\cup C_3}|I(T)|\le 1.431d_{C_2\cup C_3}\big(\log\frac{\sum_{T\in C_2\cup C_3}n(T)}{d_{C_2\cup C_3}}+1.1242\big)}
  \]
  where $d_{C_2\cup C_3}$ is the number of defective items identified by $C_2\cup C_3$.
\end{lemma}

\noindent{\bf Proof.}\, By the definition of $C_3P$, it follows that
\begin{eqnarray*}
  \sum_{T\in C_2\cup C_3}|I(T)| = \sum_{T\in C_2}|I(T)|+\sum_{T\in C_3}|I(T)| = \sum_{T\in C_2}|I(T)|+\sum_{P\in C_3P}|I(P)|.
\end{eqnarray*}
We have
\begin{eqnarray}
& & \sum_{T\in C_2}|I(T)|+\sum_{P\in C_3P}|I(P)|\nonumber\\
&=& \sum_{T\in C_2}1.431 \cdot 1 \cdot \big(\log \frac{n(T)}{1}+1.1242\big)+\sum_{P\in C_3P}1.431d_P\big(\log \frac{n(P)}{d_P}+1.1242\big)\nonumber\\
&\le& 1.431d_{C_2\cup C_3}\big(\log\frac{\sum_{T\in C_2\cup C_3}n(T)}{d_{C_2\cup C_3}}+1.1242\big)\label{C2C3-1},
\end{eqnarray}
where equation (\ref{C2C3-1}) holds by Lemma~\ref{lem:com-d}.    
Therefore, the lemma is proved.\hfill $\Box$

\

The following lemma gives a bound on the sum of the last three terms in the formula (\ref{Zualltests}).

\begin{lemma}\label{C2C3C4}
  For $C_2\cup C_3\cup C_4$, we have
  $$
  \sum_{T\in C_2\cup C_3\cup C_4}|I(T)|\le 1.431d_{Z^u}\big(\log\frac{\sum_{T\in C_2\cup C_3\cup C_4}n(T)}{d_{Z^u}}+1.1242\big)+16,
  $$
  where $3\le d_{Z^u}\le \sum_{T\in C_2\cup C_3\cup C_4}n(T)$.
\end{lemma}

\noindent{\bf Proof Sketch.}\,~The proof is technical, and the complete details are available in Appendix~\ref{sec:appendix2}. Since the tests in $C_4$ are all of the results pure, even if $^1T^1 \in C_4$, $C_4$ alone has no desired bound associated with $d_{Z^u}$. As such, we first extract a carefully designed subset $C'$ from $C_2\cup C_3$ and merge it with $C_4$ as a whole. By using this method, we bound $C'\cup C_4$ and $C_2\cup C_3 \setminus C'$ separately. Finally, obtain an upper bound on $C_2\cup C_3\cup C_4$ by Lemma~\ref{lem:com-d}. \hfill$\Box$

\subsection{Bounding the Number of Tests Performed by Strategy Up-Zig-Zag}\label{subsec:UBsZu}

{\noindent\bf Proof of Lemma~\ref{lem:UpperBound1-Zu}.} \, Based on Lemma~\ref{C_1} and Lemma~\ref{C2C3C4}, we have
\begin{eqnarray*}
  M_{Z^u}(d|n) ~=~ t_{Z^u} ~=~ \sum_{T\in \mathscr{T}}|I(T)| &=& \sum_{T\in C_2\cup C_3\cup C_4}|I(T)|+\sum_{T\in C_1}|I(T)| \\
  &\le& 1.431d_{Z^u}\big(\log\frac{\sum_{T\in \mathscr{T}}n(T)}{d_{Z^u}}+1.1242\big)+16+7 \\
  &=& 1.431d \big(\log\frac{n}{d}+1.1242\big)+23.
\end{eqnarray*}
and so Lemma~\ref{lem:UpperBound1-Zu} is proved.\hfill$\Box$

\

{\noindent\bf Proof of Lemma~\ref{lem:UpperBound2-Zu}.} \, We prove this lemma by induction on $d$ for every test class.

For $T\in C_1$, since $|I(T)| = 1$, it follows that
\[\begin{split}
\sum_{T\in C_1}|I(T)| \le \min\Big\{\sum_{T\in C_1} n(T), 7\Big\} \le \min\Big\{1.4\sum_{T\in C_1} n(T), 7\Big\}.
\end{split}\]

For $T\in C_2$, each test $T$ is applied to test out one contaminated item on a subset $S$ of size $1$. Thus, we have
\[\begin{split}
  \sum_{T\in C_2}|I(T)| = \sum_{T\in C_2} n(T) \le 1.4\sum_{T\in C_2} n(T). 
\end{split}\]

For $T\in C_3$, also $P\in C_3P$, consider the following types of $P$ as defined in Claim~\ref{ob:specialtuple}.
\begin{itemize}
  \item[(a)]~$P=(^0T^0,^2T^1,\emptyset)$, then $|I(P)| = 4$ and $n(P)=3$. Thus, $|I(P)| \le 1.4 n(P)$;
  \item[(b)]~$P=(^0T^1,^1T^2,\emptyset)$, then $|I(P)| \le 4$ and $n(P)\ge 3$. Thus, $|I(P)| \le 1.4 n(P)$;
  \item[(c)]~$P=(^1T^1,^1T^2\setminus ^2T^2\setminus ^3T^2,\emptyset)$, then $|I(P)| = 7$ and $n(P)=5$. Thus, $|I(P)| \le 1.4 n(P)$;
  \item[(d)]~$P=(T^{i-1}, T^{i},T')$ for $i\ge 3$. Let $\Delta_P = 1.4 n(P) - |I(P)|$. Assume that $T^i$ is incurred by $S^i_b$ and apply procedure 4-Split to it. If $Y$ in $S^i_b$ is contaminated, then $|I(P)|$ is no more than $i+2$ with $n(P) \ge 3\cdot 2^{i-3} + 1$. Thus,
    \begin{eqnarray*}
      \Delta_{P} \ge 1.4\cdot(3\cdot 2^{i-3} + 1) - (i+2) = 4.2\cdot 2^{i-3} - i - 0.6,
    \end{eqnarray*} 
  by elementary calculations that $\Delta_{P}\ge0$ for $i\ge3$. Analogously, if $Z$, $U$, or $V$ in $S^i_b$ is contaminated, we can get that $\Delta_P \ge 0$ (by Claim~\ref{ob:A}). So for $P\in C_3P$, $|I(P)| \le 1.4 n(P)$.
\end{itemize} 
Based on (a), (b), (c), and (d), we have $\sum_{T\in C_3}|I(T)| \le 1.4\sum_{T\in C_3} n(T)$. 

For $T\in C_4$, if $^1T^1 \notin C_4$, $|I(T)| = 1$ and $n(T)\ge1$, we have
\begin{eqnarray}
  \sum_{T\in C_4} |I(T)| \le 1.4 \sum_{T\in C_4} n(P);\label{T-C4}
\end{eqnarray}
if $^1T^1 \in C_4$, $\sum_{T\in C_4} |I(T)| = 1 + 3 + 1 + \cdots + 1 \le 1.4\times(1 + 2 + 3 +\cdots+ 3\cdot 2^{i-2}) \le 1.4\sum_{T\in C_4} n(P)$, implying that it also satisfies the formula~(\ref{T-C4}).

Therefore, by combining the above cases on $C_1,C_2,C_3$, and $C_4$, it is intuitive to get that
\begin{eqnarray*}
M_{Z^u}(d|n) \le \sum_{T\in C_1\cup C_2\cup C_3\cup C_4} |I(T)| = \sum_{T\in \mathscr{T}} n(T) \le 1.4n, 
\end{eqnarray*}
which complete the proof of Lemma~\ref{lem:UpperBound2-Zu}.\hfill$\Box$

\section{Final Remarks}\label{sec:Cr}

We propose a symmetric combinatorial algorithm for solving the competitive group testing problem. The algorithm framework primarily comprises two procedures: an original zig-zag strategy and a novel up-zig-zag strategy. Notably, the up-zig-zag strategy introduces a new testing procedure that has not been explored previously. This innovative iterative procedure enables the efficient identification of contaminated subsets, optimizing testing efficiency. On the theoretical aspect, our new algorithm does push the competitive ratio down to $1.431$. We also point out that exclusively applying the up-zig-zag strategy to the CGT-U achieves a competitive ratio of $1.431$ when $d=0$ and $3\le d< n$. Due to the straightforward implementation of the up-zig-zag strategy, we recommend it over Algorithm $Z^{c}$ for practical applications. 

Our work provides a foundation for future studies that can explore various promising directions. First, it is of great interest to improve the competitive ratio for the CGT-U. Indeed, the ratio of $1.431$ is tightly bounded by testing the number of tests for subsets of size two and three, i.e., Case 1b in the proof of Lemma~\ref{C_3}. In a future study, one can attempt to improve Procedure~\ref{strategy:S2} and/or Procedure~\ref{strategy:S3}, and thus a better approximation algorithm can be obtained. Second, it would be worthwhile to explore the potential of the up-zig-zag framework with some other novel grouping strategies to design a better strongly competitive algorithm for the CGT-U.

In addition, another promising research direction is to integrate an adaptive estimation phase for the number of defectives $d$ before applying the main testing procedure. Specifically, one could first estimate $d$ exactly or approximately using $O(d\log d)$ or $O(\log\log d)$ tests \citep{cheng2011efficient, falahatgar2016estimating}, and then apply existing algorithms designed for the case where $d$ is known to recover all defectives. This hybrid approach may potentially achieve better performance, provided that the estimation step can be efficiently derandomized. Finally, whether there exists a nontrivial lower bound other than $\lceil \log \binom{n}{d} \rceil$ for combinatorial group testing with an unknown number of defectives stands as a central open problem in this field. Investigating new heuristics with improved practical performances is also an intriguing area for further exploration.

\

\noindent \textbf{\Large {Acknowledgments}}

\

This work was supported by the China Scholarship Council (202206280182).


\renewcommand{\refname}{References}

\bibliography{ref}

\newpage

\noindent \textbf{\Large {Appendix}}

\appendix
\setcounter{table}{0}
\setcounter{figure}{0}
\setcounter{observation}{0}
\setcounter{equation}{0}
\renewcommand{\thetable}{A\arabic{table}}
\renewcommand{\thetable}{B\arabic{table}}
\renewcommand{\thefigure}{A\arabic{figure}}
\renewcommand{\theobservation}{A\arabic{observation}}
\renewcommand{\theobservation}{B\arabic{observation}}
\renewcommand{\theequation}{A.\arabic{equation}}
\renewcommand{\theequation}{B.\arabic{equation}}

\renewcommand{\thelemma}{B\arabic{lemma}}

\begin{appendix}

\section{Proof of Claim~\ref{ob:A}}\label{sec:appendix1}

Since we execute procedure 4-Split to $D(T_x^{i_x})$ and partition $T_x^{i_x}$ into four subsets $Y_x, Z_x, U_x, V_x$, so there exist the following four types of test outcomes on $P_x$.

The first type {\bf T1} is based on the case when $Y_x$ is contaminated. In this case, $I(P_x)$ contains at most four tests ($T_x^{i_x-1}$, $T_x^{i_x}$, possibly $T'_x$ and one test on $Y_x$) together with tests performed by DIG$(Y_x)$. Since $|Y_x|\le 2^{i_x-2}$, by Lemma~\ref{lem:DIG}, DIG$(Y_x)$ uses at most $i_x-2$ tests. Hence,
$
|I(P_x)| \le |I(T'_x)| + i_x+1.
$
If $T_x^{i_x-1}$ is the sixth test in its phase and $|D(T_x^{i_x})| = 3\cdot 2^{i_x-2}$, then $T'_x \neq \emptyset$ and $|I(T'_x)| = 1$; otherwise, $T'_x=\emptyset$ and $|I(T'_x)| = 0$. On the other hand, $D(T_x^{i_x-1})=3\cdot 2^{i_x-3}$, and so one defective item from $Y_x$ and $3\cdot 2^{i_x-3}$ good items from $S_a$ are identified, implying that
\begin{eqnarray*}
n(P_x) \ge |D(T_x^{i_x-1})|+1= 3\cdot2^{i_x-3}+1.
\end{eqnarray*}

The second type {\bf T2} is based on the case when $Y_x$ is pure and $Z_x$ is contaminated. In this case, $I(P_x)$ contains at most five tests ($T_x^{i_x-1}$, $T_x^{i_x}$, one test on $Y_x$, and possibly $T'_x$ and one more test on $Z_x$) together with tests performed by DIG$(Z_x)$. Since $|Z_x|\le 2^{i_x-2}$, DIG$(Z_x)$ uses at most $i_x-2$ tests. Hence, 
$
|I(P_x)| \le |I(T'_x)| + i_x+2.
$
Similar to {Case 1b}, if $T_x^{i_x-1}$ is the sixth test in the current phase and $|D(T_x^{i_x})| = 3\cdot 2^{i_x-2}$, $|I(T'_x)| = 1$; otherwise, $|I(T'_x)| = 0$. Since $Z_x$ is not empty, then $|Y_x| = 2^{i_x-2}$. So one defective item from $Z_x$, $2^{i_x-2}$ good items from $Y_x$, together with $3\cdot 2^{i_x-3}$ good items from $D(T_x^{i_x-1})$ are identified, implying that
\begin{eqnarray*}
  n(P_x) \ge |D(T_x^{i_x-1})|+|Y_x|+1=5\cdot2^{i_x-3}+1.
\end{eqnarray*}

The third type {\bf T3} is based on the case when $Y_x$ and $Z_x$ are pure, $U_x$ is contaminated. In this case, $I(P_x)$ contains at most six tests ($T_x^{i_x-1}$, $T_x^{i_x}$, two tests on $Y_x$ and $Z_x$ respectively, and possibly $T'_x$ and one more test on $U_x$) together with tests performed by DIG$(U)$. Since $|U_x|\le 2^{i_x-3}$, DIG$(U_x)$ uses at most $i_x-3$ tests. Hence, $|I(P_x)| \le |I(T'_x)| + i_x+2$. Because $U_x$ is not empty, then $|Y_x|=|Z_x|=2^{i_x-2}$. So one defective item from $U_x$, $2^{i_x-2}$ good items from $Y_x$, $2^{i_x-2}$ good items from $Z_x$, together with $3\cdot 2^{i_x-3}$ good items from $D(T_x^{i_x-1})$ are identified, implying that
\begin{eqnarray*}
  n(P) \ge |D(T_x^{i_x-1})|+|Y_x|+|Z_x|+1=7\cdot2^{i-3}+1.
\end{eqnarray*}

The fourth type {\bf T4} is based on the case when $Y_x$, $Z_x$, and $U_x$ are pure, $V_x$ is contaminated. Because $D(T_x^{i_x})$ is contaminated when procedure 4-Split is applied to it, after $Y_x, Z_x$ and $U_x$ are all tested to be pure, the last subset $V_x$ is known to be contaminated, and a test on it is not necessary. Thus, in this case, $I(P_x)$ contains at most six tests ($T_x^{i_x-1}$, $T_x^{i_x}$, three tests on $Y_x,Z_x$, and $U_x$ respectively, and possibly $T'_x$) together with tests performed by DIG$(Z)$. Since $|V_x|\le 2^{i_x-3}$, $|I(P_x)|\le |I(T'_x)|+i_x+2$. Compared with {Case 2c}, the number of tests in these two cases is the same, whereas the number of identified items in {Case 2d} is more than {Case 2c}. Since $V_x$ is not empty, then $|Y_x|=|Z_x|=2^{i_x-2}$ and $|U_x|=2^{i_x-3}$. So one defective item from $V_x$ is identified, and all items in $Y_x,Z_x,U_x$, and $D(T_x^{i_x-1})$ are identified as good items, implying that
\begin{eqnarray*}
  n(P_x) \ge |D(T_x^{i_x-1})|+|Y_x|+|Z_x|+|U_x|+1=8\cdot2^{i_x-3}+1.
\end{eqnarray*}

To sum up, Claim~\ref{ob:A} is established. Notice that compared with {\bf T2}, {\bf T3} has the same number of tests but can test more items, implying that {\bf T3} has better test results than {\bf T2}. And the same is true for {\bf T4} compared with {\bf T3}.\hfill $\Box$

\section{Proof of Lemma~\ref{C2C3C4}}\label{sec:appendix2}

Since $d_{Z^u} \ge 3$, there must exist at least three tests in $C_2\cup C_3$ or $C_2\cup C_3\cup C_4$. Construct a test subset $C' \subseteq C_2\cup C_3P \subseteq \mathscr{T}$, such that the tests in $C_4\cup C'$ can identify exactly three contaminated items. If $\sum_{T\in C_4\cup C'} |I(T)| \le 1.431 \cdot 3 \big(\log \frac{\sum_{T\in C_4\cup C'}n(T)}{3}+1.1242\big)$ + x, we can get $ \sum_{T\in C_2\cup C_3\cup C_4} |I(T)|$
  {\small
  \begin{align}
    =& \sum_{T\in C_2\cup C_3P\setminus C'} |I(T)| + \sum_{T\in C_4\cup C'} |I(T)| \nonumber\\
    \le& 1.431(d_{Z^u}-3)\big(\log \frac{\sum_{T\in C_2\cup C_3P\setminus C'}n(T)}{d_{Z^u}-3} + 1.1242\big) + 1.431 \cdot 3 \big(\log \frac{\sum_{T\in C_4\cup C'}n(T)}{3}+1.1242\big) + 16 \label{C'-0}\\
    \le& 1.431d_{Z^u}\big(\log \frac{\sum_{T\in C_2\cup C_3\cup C_4}n(T)}{d_{Z^u}}+1.1242\big)+x \nonumber,
  \end{align}}
so the Lemma~\ref{C2C3C4} holds. Since the first term of formula~(\ref{C'-0}) is supported by Lemma~\ref{C2C3}, it is only necessary to prove that the second term of formula~(\ref{C'-0}) is valid. In the following, we define $d_{C_4\cup C'}$ as the number of contaminated items identified by tests in $C_4\cup C'$, i.e., $d_{C_4\cup C'}=3$.
  
We now enumerate all the combinations of $(C_4, C')$ to verify this lemma.
  
\
  
\noindent{{\bf $\bullet$~Case 1.}}~$C_4$ is empty, i.e., $|C_4| = 0$.
  
Set $C' = \emptyset$. Based on Lemma~\ref{C2C3}, we have
  \begin{eqnarray*}
    \sum_{T\in C_2\cup C_3\cup C_4}|I(T)| &=& \sum_{T\in C_2}|I(T)|+\sum_{T\in C_3}|I(T)|\\
    &\le& 1.431d_{Z^u}\big(\log \frac{\sum_{T\in C_2\cup C_3\cup C_4}n(T)}{d_{Z^u}}+1.1242\big).
  \end{eqnarray*}
  
\
  
{\noindent{\bf $\bullet$~Case 2.}}~$|C_4| \ge 1$ and $^1T^1\notin C_4$. 
  
By Claim~\ref{ob:rank}, let $(0, 1, \cdots, i_\gamma)$ be the ranks of tests in $C_4$ for $i_\gamma \ge 0$, so $\sum_{T\in C_4} n(T) = 3\cdot 2^{i_\gamma-1}$ and $\sum_{T\in C_4}|I(T)| = i_\gamma + 1$. Since $^1T^1\notin C_4$, each $T$ in $C_4$ satisfies $|I(T)| = 1$ and $C_3P \neq \emptyset$. Assume that $C_3P$ is empty, then all of the contaminated tests are in $C_2$. Thus, $C_4$ is an empty set, which contradicts that $|C_4| \ge 1$ in this case. 
  
Let $i_m$ be the largest rank for all contaminated tests in $\mathscr{T}$ and so $i_m \ge 2$. Then we have 
  \begin{eqnarray}
  i_m = i_\gamma + \Delta_\gamma \label{equ:m=g+d},
  \end{eqnarray}
where (\ref{equ:m=g+d}) holds by Claim~\ref{ob:rank} and $\Delta_\gamma \ge 2$. Define $P_m$ as a tuple with the largest rank in $C_3P$, such that $P_m = (T^{i_m-1}_m, T^{i_m}_m, T'_m)$, in which $T^{i_m-1}_m$ has test result pure and is applied to a subset $S_m^1$ of size exactly $a_{i_m-1}$, $T^{i_m}_m$ has test result contaminated and is applied to a subset $S_m^2$ of size at most $a_{i_m}$, and $|I(T'_m)| = 0$ or $1$. Next, we analyze $i_m$. Define a function $h$ on $\sum_{T\in 
C_4\cup C'} n(T)$ and $d_{C_4\cup C'}$ such that 
  $$
  h = 1.431\cdot d_{C_4\cup C'}\big(\log \frac{\sum_{T\in 
  C_4\cup C'} n(T)}{d_{C_4\cup C'}} + 1.1242\big).
  $$
  
{\noindent{\bf -~Case 2a.}}~$i_m = 2$. 
  
By Claim~\ref{def:rank}, $C_4$ contains at most one test of rank $0$. Set $C' = \emptyset$ and we have 
  $$
  \sum_{T\in C_2\cup C_3\cup C_4}|I(T)| \le 1.4315d_{Z^u}\big(\log\frac{\sum_{T\in C_2\cup C_3\cup C_4}n(T)}{d_{Z^u}}+1.1242\big)+1.
  $$ 
  
{\noindent{\bf -~Case 2b.}}~$i_m \ge 3$. 
  
Set $C'$ to be the different sets depending on the value of $|C_1|$.
  
{{\bf $\ast$~Case 2b1.}}~$|C_1| \ge 2$. 
  
Given two arbitrary tests $T^0_\alpha,~T^0_\beta \in C_3$, i.e., $T_\alpha^0$ is applied to a subset $S_{\alpha}$ with a contaminated test result of size exactly $1$, and $T_\beta^{0}$ is applied to another subset $S_{\beta}$ with a contaminated test result of size exactly $1$. In this case, set $C'=\{T^0_\alpha, T^0_\beta, P_m\}$. Apply procedure 4-Split to $S^2_m$ and partition it into four subsets $Y_m$, $Z_m$, $U_m$, and $V_m$. There are four possible test results on $S^2_m$ by Claim~\ref{ob:A}.
  
{{\bf(\rmnum{1})} $Y_m$ is contaminated}. On the other hand, $|S^1_m| = 3\cdot 2^{i_m-3}$, and so three defective items from $Y_m,S_\alpha$ and $S_\beta$ respectively are identified, and all items tested by $C_4$ and all items in $S^1_m$ are identified as good items, implying that
  \begin{eqnarray*}
    \sum_{T\in C_4\cup C'} n(T) &=& \sum_{T\in C_4} n(T) + n(P_m) + n(T_\alpha^{0}) + n(T_\beta^{0})\\
    &\ge& \sum_{T\in C_4} n(T) + |S^1_m| + 1 + n(T_\alpha^{0}) + n(T_\beta^{0})\\
    &=& 3\cdot 2^{i_\gamma-1} + 3\cdot 2^{i_m-3} + 3.
  \end{eqnarray*}
Thus, $h$ satisfies that
  \begin{eqnarray}
    h &\ge& 1.431\cdot 3\Big(\log\frac{3\cdot 2^{i_\gamma-1}+3\cdot 2^{i_m-3}+3}{3}+1.1242\Big)\nonumber\\
    &=& 4.293\cdot \Big(\log\big(\frac{1}{2^{\Delta_\gamma+1}}+\frac{1}{2^{i_m}}+\frac{1}{8}\big)+i_m+1.1242\Big)\label{2-2-1-0}\\
    &\ge& 4.293i_m-8.1, \label{2-2-1-0-1}
  \end{eqnarray}
where (\ref{2-2-1-0}) holds by (\ref{equ:m=g+d}) and (\ref{2-2-1-0-1}) holds by $\Delta_\gamma \ge 2$.
  
In this case, $\bigcup_{T\in C_4\cup C'}I(T)$ contains at most six tests ($T^{i_m-1}_m$, $T^{i_m}_m$, $T^0_\alpha$, $T^0_\beta$; possibly $T'_m$ and one test on $Y_m$) together with tests in $C_4$ and performed by DIG($Y_m$). Since $|Y_m|\le 2^{i_m-2}$, DIG($Y_m$) uses at most $i_m-2$ tests. So the total number of tests in $C_4\cup C'$ satisfies that
  \begin{eqnarray*}
    \sum_{T\in C_4\cup C'}|I(T)| &=& \sum_{T\in C_4}|I(T)| + |I(P_m)| + |I(T_\alpha^{0}) + |I(T_\beta^{0})| \\
    &\le& (i_\gamma + 1) + (|I(T'_m)|+ i_m + 1) + 2 \\
    &=& 2i_m - \Delta_\gamma + 4 + |I(T'_m)|.
  \end{eqnarray*}
Thus,
  \begin{eqnarray}
  h-\sum_{T\in C_4\cup C'}|I(T)| &\ge& (4.293i_m-8.1)-(2i_m - \Delta_\gamma + 4 + |I(T'_m)|)\nonumber\\ 
  &\ge& 2.293i_m+\Delta_\gamma-|I(T'_m)|-12.1.\label{incre1}
  \end{eqnarray}
Based on $\Delta_\gamma\ge2$, if $3\le i_m < 6$, $|I(T'_m)| = 0$ and so $h - \sum_{T\in C_4\cup C'}|I(T)| \ge -4$; if $i_m \ge 6$, $|I(T'_m)| \le 1$ and so $h-\sum_{T\in C_4\cup C'}|I(T)| > 2$. Thus, $h-\sum_{T\in C_4\cup C'}|I(T)| > -4$ for $i_m\ge 3$.
  
{{\bf(\rmnum{2})} $Y_m$ is pure and $Z_m$ is contaminated}. So three defective items from $Z_m,S_\alpha$, and $S_\beta$ respectively are identified, and all items tested by $C_4$ and all items in $S^1_m, Y_m$ are identified as good items, implying that
  \begin{eqnarray*}
    \sum_{T\in C_4\cup C'}n(T) &\ge& 3\cdot 2^{i_\gamma-1} + 5\cdot 2^{i_m-3} + 3.
  \end{eqnarray*}
  
In this case, $\bigcup_{T\in C_4\cup C'}I(T)$ contains at most five tests ($T^{i_m-1}_m$, $T^{i_m}_m$, one test on $Y_m$; possibly $T'_m$ and one more test on $Z_m$) together with tests in $C_4$ and performed by DIG($Z_m$). Since $|Z_m|\le 2^{i_m-2}$, DIG($Z_m$) uses at most $i_m-2$ tests. Similar to case {\bf(\rmnum{1})}, we have
  \begin{eqnarray*}
    \sum_{T\in C_4\cup C'}|I(T)| &\le& 2i_m - \Delta_\gamma + 5 + |I(T'_m)|,
  \end{eqnarray*}
and
  \begin{eqnarray*}
    h \ge 1.431\cdot 3\cdot \Big(\log\frac{3\cdot 2^{i_\gamma-1} + 5\cdot 2^{i_m-3}+3}{3} + 1.1242\Big) > 4.293i_m - 4.9.
  \end{eqnarray*}
Thus,
  \begin{eqnarray}
    h-\sum_{T\in C_4\cup C'}|I(T)| 
    &>& 2.293i_m + \Delta_\gamma - |I(T'_m)| - 9.9.\label{incre2}
  \end{eqnarray}
    
By comparing inequalities (\ref{incre1}) and (\ref{incre2}), it is easy to get that the increment of $h$ is greater than the increment of $\sum_{T\in C_4\cup C'} |I(T)|$ when the test result of $P_m$ from {\bf T1} changes to {\bf T2} according to Claim~\ref{ob:A}. Thus, $h - \sum_{T\in C_4\cup C'} |I(T)| > -4$ in this case. Similar to the analysis of the test result that $Y_m$ or $Z_m$ is contaminated, ‘$h-\sum_{T\in C_4\cup C'}|I(T)|$’ is definitely larger than ‘$-4$’ when $U$ or $V$ is contaminated.
  
Hence in {Case 2b1}, we have
  \begin{eqnarray}
    \sum_{T\in C_4 \cup C'} |I(T)| &\le& 1.431 d_{C_4\cup C'}\big(\log \frac{\sum_{T\in C_4 \cup C'} n(T)}{d_{C_4\cup C'}} + 1.1242\big)+ 4.\label{2-2-1-3}
  \end{eqnarray}
  
  {{\bf $\ast$~Case 2b2.}}~$|C_2|=1$. 
  
Given $T_\alpha^{0} \in C_2$, i.e., $T_\alpha^{0}$ is applied to a subset $S_\alpha$ of size exactly $1$. Given a tuple $P_\beta \in C_3P$, i.e., $P_\beta = (T_\beta^{i_\beta-1}, T_\beta^{i_\beta}, T'_\beta)$, in which $T_\beta^{i_\beta-1}$ has a pure test result and is applied to a subset $S^1_\beta$ of size exactly $a_{i_\beta-1}$, $T_\beta^{i_\beta}$ has test result contaminated and is applied to a subset $S^2_\beta$ of size at most $a_{i_\beta}$, and $|I(T'_\beta)| = 0$ or $1$. Since $i_m$ is the largest rank in $C_3P$, we have that 
  \begin{eqnarray}
    i_m = i_\beta+\Delta_\beta, \label{equ:m=b+D} 
  \end{eqnarray}
where $\Delta_\beta \ge 0$. In this case, set $C'=\{T_\alpha^{0},P_\beta,P_m\}$. Consider two possible cases below on $i_\beta$.
  
{If $1 \le i_\beta \le 2$}. By Claim~\ref{ob:specialtuple}, it only increases no more than six tests incurred by $P_\beta$ compared with $T_\beta^{0}$ in {Case~2b1}. Thus, based on (\ref{2-2-1-3}), we can get
  \begin{eqnarray*}
    \sum_{T\in C_4\cup C'} |I(T)| \le 1.431 d_{C_4\cup C'}\big(\log \frac{\sum_{T\in C_4\cup C'} n(T)}{d_{C_4\cup C'}} + 1.1242\big) + 10.
  \end{eqnarray*}
  
{If $i_\beta \ge 3$}. Apply procedure 4-Split to $S^2_\beta$ and partition it into four subsets $Y_\beta,Z_\beta,U_\beta$ and $V_\beta$. So there are four possible test results {\bf T1}-{\bf T4} on $P_\beta$ based on Claim~\ref{ob:A}.
  
{{\bf(\rmnum{1})}~$Y_\beta$ is contaminated.}
     
{{\bf(\rmnum{1}-1)}~$Y_m$ is contaminated}. So three defective items from $Y_m$, $Y_\beta$, and $S_\alpha$ respectively are identified, and all items tested by $C_4$, all items in $S_m^1$ and in $S_\beta^1$ are identified as good items, implying that
    \begin{eqnarray*}
      \sum_{T\in C_4\cup C'} |n(T)| &=& \sum_{T\in C_4}n(T) + n(P_m) + n(P_\beta) + n(T_\alpha^{0})\\
      &\ge& \sum_{T\in C_4}n(T) + |S^1_m|+1+|S^1_\beta|+1+n(T^0_\alpha)\\
      &=& 3\cdot 2^{i_\gamma-1} + 3\cdot 2^{i_m-3} + 3\cdot 2^{i_\beta - 3} + 3.
    \end{eqnarray*}
In this case, $\bigcup_{T\in C_4\cup C'}I(P)$ contains at most eight tests ($T^{i_m-1}_m$, $T^{i_m}_m$, $T^{i_\beta-1}_\beta$, $T^{i_\beta}_\beta$; possibly $T'_m$, $T'_\beta$, and two more tests on $Y_m, Y_\beta$) together with tests in $C_4$, and performed by DIG($Y_m$) and DIG($Y_\beta$). Since $|Y_m|\le 2^{i_m-2}$ and $|Y_\beta| \le 2^{i_\beta-2}$, DIG($Y_m$) and DIG($Y_\beta$) use at most $i_m-2$ and $i_\beta-2$ tests, respectively. So we have
    \begin{eqnarray*}
      \sum_{T\in C_4\cup C'} |I(T)| &=& \sum_{T\in C_4}|I(T)| + |I(P_m)| + |I(P_\beta)| + |I(T_\alpha^{0})|\\
      &\le& (i_\gamma + 1) + (i_m + 1 + |I(T'_m)|) + (i_\beta + 1 + |I(T'_\beta)|) + 1\\
      &=& 3i_m - \Delta_\gamma - \Delta_\beta + |I(T'_m)| + |I(T'_\beta)| + 4.
    \end{eqnarray*}
Thus,
    \begin{eqnarray}
      h &\ge& 1.431\cdot 3\Big(\log\frac{3\cdot 2^{i_\gamma-1} + 3\cdot 2^{i_m-3} + 3\cdot 2^{i_\beta - 3} + 3}{3}+1.1242\Big)\nonumber\\
      &=& 4.293\cdot \big(\log\big(\frac{1}{2^{\Delta_\gamma+1}} + \frac{1}{2^{i_m}}+\frac{1}{2^{\Delta_\beta + 3}} + \frac{1}{8}\big)+i_m + 1.1242\Big)\label{2-2-i-i1-1}\\
      &>& 4.293i_m - 8.1\label{2-2-i-i1-2},
    \end{eqnarray}
where (\ref{2-2-i-i1-1}) holds by (\ref{equ:m=g+d}) and (\ref{equ:m=b+D}), and (\ref{2-2-i-i1-2}) holds by $\Delta_\gamma \ge 2$ and $\Delta_\beta \ge 0$. For $i_m\ge 3$, we have
    \begin{eqnarray}
      h - \sum_{T\in C_4\cup C'} |I(T)| &\ge& (4.293i_m - 8.1) - (3i_m - \Delta_\gamma - \Delta_\beta + |I(T'_m)| + |I(T'_\beta)| + 4)\nonumber \\
      &=& 1.293i_m + \Delta_\gamma + \Delta_\beta - |I(T'_m)| - |I(T'_\beta)| - 12.1 \label{incre3} \\
      &>& -9.\nonumber
    \end{eqnarray}
  
{{\bf(\rmnum{1}-2)}~$Y_m$ is pure and $Z_m$ is contaminated}. So three defective items from $Z_m$, $Y_\beta$, and $S_\alpha$ respectively are identified, and all items tested by $C_4$, all items in $S_m^1$, in $S_\beta^1$, and in $Y_m$ are identified as good items, implying that
    \begin{eqnarray*}
      \sum_{T\in C_4\cup C'} |n(T)| &\ge& \sum_{T\in C_4}n(T) + |Y_m| + |S^1_m| + 1 + |S^1_\beta| + 1 + n(T_\alpha^{0})\\
      &\ge& 3\cdot 2^{i_\gamma-1} + 5\cdot 2^{i_m-3} + 3\cdot 2^{i_\beta - 3} + 3.
    \end{eqnarray*}
In this case, $\bigcup_{T\in C_4\cup C'}I(P)$ contains at most nine tests ($T^{i_m-1}_m$, $T^{i_m}_m$, $T^{i_\beta-1}_\beta$, $T^{i_\beta}_\beta$, one test on $Y_m$; possibly $T'_m$, $T'_\beta$, and two more tests on $Y_\beta, Z_\beta$) together with tests in $C_4$, and performed by DIG($Y_m$) and DIG($Z_\beta$). Since $|Z_m|\le 2^{i_m-2}$ and $|Y_\beta| \le 2^{i_\beta-2}$, DIG($Z_m$) and DIG($Y_\beta$) use at most $i_m-2$ and $i_\beta-2$ tests, respectively. So we have
    \begin{eqnarray*}
      \sum_{T\in C_4\cup C'} |I(T)| &\le& (i_\gamma + 1) + (i_m + 2 + |I(T'_m)|) + (i_\beta + 1 + |I(T'_\beta)|) + 1\\
      &\le& 3i_m - \Delta_\gamma - \Delta_\beta + |I(T'_m)| + |I(T'_\beta)| + 5.
    \end{eqnarray*}
Thus,
    \begin{eqnarray*}
      h \ge 1.431\cdot 3\cdot\Big(\log \frac{3\cdot 2^{i_\gamma-1} + 5\cdot 2^{i_m-3} + 3\cdot 2^{i_\beta - 3} + 3}{3} + 1.1242\Big) > 4.293i_m - 4.9.
    \end{eqnarray*}
So we have
    \begin{eqnarray}
      h - \sum_{T\in C_4\cup C'} |I(T)| &>& 1.293i_m + \Delta_\gamma + \Delta_\beta - |I(T'_m)| - |I(T'_\beta)| - 9.9\label{incre4}
   \end{eqnarray}
      
By comparing inequalities~(\ref{incre4}) and (\ref{incre3}), we can get that the increment of $h$ is more than the increment of $\sum_{T\in C_4\cup C'} |I(T)|$ when the test result of $P_m$ from {\bf T1} changes to {\bf T2}. Thus, $h - \sum_{T\in C_4\cup C'} I(T) > -9$ in this case. Similar to the analysis of the test result that $Z_m$ is contaminated, ‘$h-\sum_{T\in C_4\cup C'} |I(T)|$’ is definitely larger than ‘$-9$’ when $U_m$ or $V_m$ is contaminated.
  
{{\bf(\rmnum{2})}~$Y_\beta$ is pure and $Z_\beta$ is contaminated.}
    
{{\bf(\rmnum{2}-1)}~$Y_m$ in $S^2_m$ is contaminated}. So three defective items from $Y_m$, $Z_\beta$, and $S_\alpha$ respectively are identified, and all items tested by $C_4$, all items in $S^1_m$, in $S^1_\beta$ and in $Y_\beta$ are identified as good items, implying that
    \begin{eqnarray*}
      \sum_{T\in C_4\cup C'} |n(T)| &\ge& \sum_{T\in C_4}n(T) + |S^1_m| + 1 + |Y_\beta| + |S^1_\beta| + 1 + n(T_\alpha^{0})\\
      &\ge& 3\cdot 2^{i_\gamma-1} + 3\cdot 2^{i_m-3} + 5\cdot 2^{i_\beta - 3} + 3.
    \end{eqnarray*}
In this case, $\bigcup_{T\in C_4\cup C'}I(P)$ contains at most nine tests ($T^{i_m-1}_m$, $T^{i_m}_m$, $T^{i_\beta-1}_\beta$, $T^{i_\beta}_\beta$, one test on $Y_\beta$; possibly $T'_m$, $T'_\beta$, and two more tests on $Y_m, Z_\beta$) together with tests in $C_4$, and performed by DIG($Y_m$) and DIG($Z_\beta$). Since $|Y_m|\le 2^{i_m-2}$ and $|Z_\beta| \le 2^{i_\beta-2}$, DIG($Y_m$) and DIG($Z_\beta$) use at most $i_m-2$ and $i_\beta-2$ tests, respectively. So we have
    \begin{eqnarray*}
      \sum_{T\in C_4cup C'} |I(T)| &\le& (i_\gamma + 1) + (i_m + 1 + |I(T'_m)|) + (i_\beta + 2 + |I(T'_\beta)|) + 1\\
      &=& 3i_m - \Delta_\gamma - \Delta_\beta + |I(T'_m)| + |I(T'_\beta)| + 5.
    \end{eqnarray*}
Thus,
    \begin{eqnarray*}
      h \ge 1.431\cdot 3\cdot\Big(\log \frac{3\cdot 2^{i_\gamma-1} + 3\cdot 2^{i_m-3} + 5\cdot 2^{i_\beta - 3} + 3}{3} + 1.1242\Big) > 4.293i_m - 8.1.\nonumber
    \end{eqnarray*}
For $i_m \ge 3$, we have
    \begin{eqnarray}
      h - \sum_{T\in C_4\cup C'} |I(T)| &>& (4.293i_m - 8.1) - (3i_m - \Delta_\gamma - \Delta_\beta + |I(T'_m)| + |I(T'_\beta)| + 5)\nonumber\\
      &=& 1.293i_m + \Delta_\gamma + \Delta_\beta - |I(T'_m)| - |I(T'_\beta)| - 13.1\label{incre5}\\
      &>& -10\nonumber.
    \end{eqnarray}
  
{{\bf(\rmnum{2}-2)}~$Y_m$ is pure and $Z_m$ is contaminated}. So three defective items from $Z_m$, $Z_\beta$, and $S_\alpha$ respectively are identified, and all items tested by $C_4$, all items in $S^1_m$, in $S^1_\beta$, in $Y_m$ and in $Y_\beta$ are identified as good items, implying that
    \begin{eqnarray*}
      \sum_{T\in C_4\cup C'} |n(T)| &\ge& \sum_{T\in C_4}n(T) + |Y_m| + |S^1_m| + 1 + |Y_\beta| + |S^1_\beta| + 1 + n(T_\alpha^{0})\\
      &=& 3\cdot 2^{i_\gamma-1} + 5\cdot 2^{i_m-3} + 5\cdot 2^{i_\beta - 3} + 3.
    \end{eqnarray*}
In this case, $\bigcup_{T\in C_4\cup C'}I(P)$ contains at most ten tests ($T^{i_m-1}_m$, $T^{i_m}_m$, $T^{i_\beta-1}_\beta$, $T^{i_\beta}_\beta$, two tests on $Y_m, Y_\beta$; possibly $T'_m$, $T'_\beta$, and two more tests on $Y_m, Z_\beta$) together with tests in $C_4$, and performed by DIG($Y_m$) and DIG($Z_\beta$). Since $|Y_m|\le 2^{i_m-2}$ and $|Z_\beta| \le 2^{i_\beta-2}$, DIG($Y_m$) and DIG($Z_\beta$) use at most $i_m-2$ and $i_\beta-2$ tests, respectively. So we have
    \begin{eqnarray*}
      \sum_{T\in C_4\cup C'} |I(T)| &\le& (i_\gamma + 1) + (i_m + 2 + |I(T'_m)|) + (i_\beta + 2 + |I(T'_\beta)|) + 1\\
      &=& 3i_m - \Delta_\gamma - \Delta_\beta + |I(T'_m)| + |I(T'_\beta)| + 6.
    \end{eqnarray*}
Thus,
    \begin{eqnarray*}
      h \ge 1.431\cdot 3\cdot\Big(\log \frac{3\cdot 2^{i_\gamma-1} + 5\cdot 2^{i_m-3} + 5\cdot 2^{i_\beta - 3} + 3}{3} + 1.1242\Big) > 4.293i_m - 4.9
    \end{eqnarray*}
So we have
      \begin{eqnarray}
        h - \sum_{T\in C_4\cup C'} |I(T)| &>&  1.293i_m + \Delta_\gamma + \Delta_\beta - |I(T'_m)| - |I(T'_\beta)| - 9.9\label{incre6}
      \end{eqnarray}
By comparing inequalities~(\ref{incre6}) and (\ref{incre5}), we can get that the increment of $h$ is more than the increment of $\sum_{T\in C_4\cup C'} |I(T)|$ when the test result of $P_m$ from {\bf T1} changes to {\bf T2}. Thus, $h - \sum_{T\in C_4\cup C'} I(T) > -10$ in this case. Similar to the analysis of the test result that $Z_m$ is contaminated, ‘$h-\sum_{T\in C_4\cup C'} |I(T)|$’ is definitely bigger than ‘$-10$’ when $U_m$ or $V_m$ is contaminated.
  
{{\bf(\rmnum{3})} $Y_\beta$ and $Z_\beta$ are pure, $U_\beta$ is contaminated.}
      
{{\bf(\rmnum{3}-1)}~$Y_m$ is contaminated}. So three defective items from $Y_m$, $U_\beta$, and $S_\alpha$ respectively are identified, and all items tested by $C_4$, all items in $S^1_m$, in $S^1_\beta$, in $Y_\beta$ and in $Z_\beta$ are identified as good items, implying that
    \begin{eqnarray*}
      \sum_{T\in C_4\cup C'} |n(T)| &\ge& \sum_{T\in C_4}n(T) + |S^1_m| + 1 + |Y_\beta| + |Z_\beta| + |S^1_\beta| + 1 + n(T_\alpha^{0})\\
      &=& 3\cdot 2^{i_\gamma-1} + 3\cdot 2^{i_m-3} + 7\cdot 2^{i_\beta - 3} + 3.
    \end{eqnarray*}
In this case, $\bigcup_{T\in C_4\cup C'}I(P)$ contains at most ten tests ($T^{i_m-1}_m$, $T^{i_m}_m$, $T^{i_\beta-1}_\beta$, $T^{i_\beta}_\beta$, two tests on $Y_\beta, Z_\beta$; possibly $T'_m$, $T'_\beta$, and two more tests on $Y_m, U_\beta$) together with tests in $C_4$, and performed by DIG($Y_m$) and DIG($U_\beta$). Since $|Y_m|\le 2^{i_m-2}$ and $|U_\beta| \le 2^{i_\beta-3}$, DIG($Y_m$) and DIG($U_\beta$) use at most $i_m-2$ and $i_\beta-3$ tests, respectively. So we have
    \begin{eqnarray*}
      \sum_{T\in C_4\cup C'} |I(T)| &\le& (i_\gamma + 1) + (i_m + 1 + |I(T'_m)|) + (i_\beta + 2 + |I(T'_\beta)|) + 1\\
      &=& 3i_m - \Delta_\gamma - \Delta_\beta + |I(T'_m)| + |I(T'_\beta)| + 5.
    \end{eqnarray*}
Thus,
    \begin{eqnarray*}
      h \ge 1.431\cdot 3\cdot\Big(\log \frac{3\cdot 2^{i_\gamma-1} + 3\cdot 2^{i_m-3} + 7\cdot 2^{i_\beta - 3} + 3}{3} + 1.1242\Big) > 4.293i_m - 8.1.\nonumber
    \end{eqnarray*}
For $i_m \ge 3$, we have
    \begin{eqnarray}
      h - \sum_{T\in C_4\cup C'} |I(T)| &\ge& (4.293i_m - 8.1) - (3i_m - \Delta_\gamma - \Delta_\beta + |I(T'_m)| + |I(T'_\beta)| + 5)\nonumber\\
      &=& 1.293i_m + \Delta_\gamma + \Delta_\beta - |I(T'_m)| - |I(T'_\beta)| - 13.1\label{incre7}\\
      &>& -10.\nonumber
    \end{eqnarray}
  
{{\bf(\rmnum{3}-2)}~$Y_m$ is pure and $Z_m$ is contaminated}. So three defective items from $Z_m$, $U_\beta$, and $S_\alpha$ respectively are identified, and all items tested by $C_4$, all items in $S^1_m$, in $S^1_\beta$, in $Y_m$, in $Y_\beta$ and $Z_\beta$ are identified as good items, implying that
    \begin{eqnarray*}
      \sum_{T\in C_4\cup C'} |n(T)| &\ge& \sum_{T\in C_4}n(T) + |S^1_m| + 1 + |Y_\beta| + |Z_\beta| + |S^1_\beta| + 1 + n(T_\alpha^{0})\\
      &=& 3\cdot 2^{i_\gamma-1} + 5\cdot 2^{i_m-3} + 7\cdot 2^{i_\beta - 3} + 3.
    \end{eqnarray*}
In this case, $\bigcup_{T\in C_4\cup C'}I(P)$ contains at most eleven tests ($T^{i_m-1}_m$, $T^{i_m}_m$, $T^{i_\beta-1}_\beta$, $T^{i_\beta}_\beta$, three tests on $Y_m, Y_\beta$ and $Z_\beta$; possibly $T'_m$, $T'_\beta$, and two more tests on $Z_m, U_\beta$) together with tests in $C_4$, and performed by DIG($Z_m$) and DIG($U_\beta$). Since $|Z_m|\le 2^{i_m-2}$ and $|U_\beta| \le 2^{i_\beta-3}$, DIG($Z_m$) and DIG($U_\beta$) use at most $i_m-2$ and $i_\beta-3$ tests, respectively. So we have
    \begin{eqnarray*}
      \sum_{T\in C_4\cup C'} |I(T)| &\le& (i_\gamma + 1) + (i_m + 2 + |I(T'_m)|) + (i_\beta + 2 + |I(T'_\beta)|) + 1\\
      &\le& 3i_m - \Delta_\gamma - \Delta_\beta + |I(T'_m)| + |I(T'_\beta)| + 6.
    \end{eqnarray*}
Thus,
    \begin{eqnarray*}
      h \ge 1.431\cdot 3\cdot\Big(\log \frac{3\cdot 2^{i_\gamma-1} + 5\cdot 2^{i_m-3} + 7\cdot 2^{i_\beta - 3} + 3}{3} + 1.1242\Big) > 4.293i_m - 4.9
    \end{eqnarray*}
So we have
    \begin{eqnarray}
      h - \sum_{T\in C_4\cup C'} |I(T)| &>& 1.293i_m + \Delta_\gamma + \Delta_\beta - |I(T'_m)| - |I(T'_\beta)| - 10.9\label{incre8}
    \end{eqnarray}
By comparing inequalities~(\ref{incre8}) and (\ref{incre7}), it is easy to get that the increment of $h$ is more than the increment of $\sum_{T\in C_4\cup C'} |I(T)|$ when the test result of $P_m$ from {\bf T1} change to {\bf T2}. Thus, $h - \sum_{T\in C_4\cup C'} I(T) > -10$. Similar to the analysis of the test result that $Z_m$ is contaminated, ‘$h-\sum_{T\in C_4\cup C'} |I(T)|$’ is definitely larger than ‘$-10$’ when $U_m$ or $V_m$ is contaminated.
  
{{\bf(\rmnum{4})} $Y_\beta, Z_\beta$, and $U_\beta$ are pure, $V_\beta$ is contaminated.} 
    
Because after $Y_\beta, Z_\beta$, and $U_\beta$, are all tested to be pure, the last subset $V_\beta$ is known to be contaminated, and a test on it is unnecessary. Thus, the result is better than the result when $U_\beta$ is contaminated with the same test number but more identified items. Therefore, $h-\sum_{T\in C_4\cup C'} |I(T)| > -10$. Thus, for $i_\beta \ge 3$, we have
  $$
    \sum_{T\in C_4\cup C'} |I(T)| \le 1.431 d_{C_4\cup C'}\big(\log \frac{\sum_{T\in C_4 \cup C'} n(T)}{d_{C_4\cup C'}} + 1.1242\big) + 10.
  $$
Hence in {Case~2b2}, for $T\in C_4\cup C'$ we have
  \begin{eqnarray}
    \sum_{T\in C_4\cup C'} |I(T)| &\le& 1.431 d_{C_4\cup C'}\big(\log \frac{\sum_{T\in C_4 \cup C'} n(T)}{d_{C_4\cup C'}} + 1.1242\big)+ 10.\label{2-2-2-7}
  \end{eqnarray}
  
\
  
{{\bf $\ast$~Case 2b3.}}~$|C_2|=0$. 
  
There must exist another tuple in $C_3P$, here called $P_\alpha$, in addition to $P_\beta$ and $P_m$. Let $P_\alpha = (T_\alpha^{i_\alpha-1}, T_\alpha^{i_\alpha}, T'_\alpha)$, in which $T_\alpha^{i_\alpha-1}$ has test result pure and is applied to a subset $S^1_\alpha$ of size exactly $a_{i_\alpha - 1}$, $T^{i_\alpha}_\alpha$ has result contaminated and is applied to a subset $S^2_\alpha$ of size at most $a_{i_\alpha}$, and $|I(T'_\alpha)|=0$ or $1$. Similar to $i_\beta$, we have that
  \begin{eqnarray}
    i_m = i_\alpha + \Delta_\alpha,\label{equ:m=a+D}
  \end{eqnarray}
where $\Delta_\alpha \ge 0$. In this case, set $C' = \{P_\alpha, P_\beta, P_m\}$. Consider the two cases below based on $i_\alpha$.
  
{If $i_\alpha \le 2$}. Compared with the result of {Case 2b2}, $C'$ increases by no more than six tests according to Claim~\ref{ob:specialtuple}. Thus, based on (\ref{2-2-2-7}), we have
  $$
  \sum_{T\in C_4\cup C'} |I(T)| \le 1.431 d_{C_4\cup C'}\big(\log \frac{\sum_{T\in C_4\cup C'} n(T)}{d_{C_4\cup C'}} + 1.1242\big) + 16.
  $$
  
{If $i_\alpha \ge 3$}. There are four possible test results for $S^2_\alpha$ like $S^2_\beta$ from {\bf T1} to {\bf T4} by applying procedure 4-Split to it. Here, we partition $S^2_\alpha$ into four subsets $Y_\alpha,Z_\alpha,U_\alpha$, and $V_\alpha$. By comparing the value of ‘$h-\sum_{T\in C_4\cup C'} |I(T)|$’, the worst case on $P_m$ and $P_\beta$ is ({\rmnum{3}}-2) in {Case~2b2} that both $U_\beta$ and $Z_m$ are contaminated according to the analysis of {Case~2b2}. All subsequent analyses in the following are based on this fact, and so we let both $U_\beta$ and $Z_m$ be contaminated in this case. 
  
{{\bf(\rmnum{1})} $Y_\alpha$ is contaminated.}
     
So three defective items from $Y_\alpha$, $Z_m$, and $U_\beta$ respectively are identified, and all items tested by $C_4$, all items in $S^1_\alpha$, in $S^1_m$, in $S^1_\beta$, in $Y_m$, in $Y_\beta$, and in $Z_\beta$ are identified as good items, implying that
    \begin{eqnarray*}
      \sum_{T\in C_4\cup C'} |n(T)| &=& \sum_{T\in C_4}n(T) + n(P_m) + n(P_\alpha) + n(P_\beta)\\
      &\ge& \sum_{T\in C_4}n(T) + |Y_m| + |S^1_m| + 1 + |Y_\beta| + |Z_\beta| + |S^1_\beta| + 1 + |S^1_\alpha| + n(T_\alpha^{0})\\
      &=& 3\cdot 2^{i_\gamma-1} + 5\cdot 2^{i_m-3} + 3\cdot 2^{i_\alpha - 3} +7\cdot 2^{i_\beta - 3} + 3.
    \end{eqnarray*}
$I(P_\alpha)$ contains at most four tests ($T^{i_\alpha - 1}_\alpha$, $T^{i_\alpha}_\alpha$, possibly $T'_\alpha$ and one test on $Y_\alpha$) together with tests performed by DIG($Y_\alpha$). Since $|Y_\alpha| \le 2^{i_\alpha -2}$, DIG($Y_\alpha$) uses at most $i_\alpha-2$ tests. Thus, we have
    \begin{eqnarray*}
      \sum_{T\in C_4\cup C'} |I(T)| &=& \sum_{T\in C_4}|I(T)| + |I(P_m)| + |I(P_\beta)| + |I(P_\alpha)|\\
      &\le& (i_\gamma + 1) + (i_m + 2 + |I(T'_m)|) + (i_\beta + 2 + |I(T'_\beta)|)\\
      & & + (i_\alpha + 1 + |I(T'_\alpha)|)\\
      &=& 4i_m - \Delta_\alpha - \Delta_\beta - \Delta_\gamma + |I(T'_m)| + |I(T'_\alpha)| + |I(T'_\beta)| + 6,
    \end{eqnarray*}
Thus,
    \begin{eqnarray}
      h &\ge& 1.431\cdot 3\Big(\log \frac{3\cdot 2^{i_\gamma-1} + 5\cdot 2^{i_m-3} + 3\cdot 2^{i_\alpha - 3} +7\cdot 2^{i_\beta - 3} + 3}{3} + 1.1242\Big)\nonumber\\
      &=& 4.293\cdot\Big(\log\big( \frac{1}{2^{\Delta_\gamma + 1}} + \frac{1}{2^{i_m}} + \frac{1}{2^{\Delta_\alpha+3}} + \frac{7}{3\cdot 2^{\Delta_\beta + 3}} + \frac{5}{24} \big) + i_m + 1.1242\Big)\label{2-3-1-1}\\
      &>& 4.293i_m - 4.9,\label{2-3-1-2}
    \end{eqnarray}
where (\ref{2-3-1-1}) holds by (\ref{equ:m=g+d}), (\ref{equ:m=b+D}), and (\ref{equ:m=a+D}); and (\ref{2-3-1-2}) holds by $\Delta_\gamma \ge 2$, $\Delta_\beta \ge 0$, and $\Delta_\alpha \ge 0$. Then for $i_m \ge 3$,
    \begin{eqnarray*}
      h - \sum_{T\in C_4\cup C'} |I(T)| &\ge& (4.293i_m - 4.9)\\
      & & - (4i_m - \Delta_\alpha - \Delta_\beta - \Delta_\gamma + |I(T'_m)| + |I(T'_\alpha)| + |I(T'_\beta)| + 6)\\
      &>& 0.293i_m - |I(T'_m)| - |I(T'_\beta)| - |I(T'_\alpha)| - 10.9\\
      &>& -11.
    \end{eqnarray*}
  
{{\bf(\rmnum{2})} $Y_\alpha$ is pure, $Z_\alpha$ is contaminated.}
     
So three defective items from $Z_\alpha$, $Z_m$, and $U_\beta$ respectively are identified, and all items tested by $C_4$, all items in $S^1_\alpha$, in $S^1_m$, in $S^1_\beta$, in $Y_m$, in $Y_\beta$, in $Z_\beta$, and in $Y_\alpha$ are identified as good items, implying that
    \begin{eqnarray*}
      \sum_{T\in C_4\cup C'} |n(T)| &\ge& \sum_{T\in C_4}n(T) + |Y_m| + |S^1_m| + 1 + |Y_\beta| + |Z_\beta| + |S^1_\beta| + 1 \\
      & & + |Y_\alpha| + |S^1_\alpha| + n(T_\alpha^{0})\\
      &\ge& 3\cdot 2^{i_\gamma-1} + 5\cdot 2^{i_m-3} + 5\cdot 2^{i_\beta - 3} +7\cdot 2^{i_\beta - 3} + 3.
    \end{eqnarray*}
$I(P_\alpha)$ contains at most five tests ($T^{i_\alpha - 1}_\alpha$, $T^{i_\alpha}_\alpha$, one test on $Y_\alpha$; possibly $T'_\alpha$ and one test on $Z_\alpha$) together with tests performed by DIG($Z_\alpha$). Since $|Z_\alpha| \le 2^{i_\alpha -2}$, DIG($Z_\alpha$) uses at most $i_\alpha-2$ tests. Thus, we have
    \begin{eqnarray*}
      \sum_{T\in C_4\cup C'} |I(T)| &\le& (i_\gamma + 1) + (i_m + 2 + |I(T'_m)|) + (i_\beta + 2 + |I(T'_\beta)|)\\
      & & + (i_\alpha + 2 + |I(T'_\alpha)|) \\
      &=& 4i_m - \Delta_\alpha - \Delta_\beta - \Delta_\gamma + |I(T'_m)| + |I(T'_\beta)| + |I(T'_\alpha)| + 7,
    \end{eqnarray*}
Thus, for $i_m \ge 3$,
    \begin{eqnarray*}
      h - \sum_{T\in C_4\cup C'}|I(T)| &\ge& (4.293i_m - 4.9)\\
      & &  - (4i_m - \Delta_\alpha - \Delta_\beta - \Delta_\gamma + |I(T'_m)| + |I(T'_\beta)| + |I(T'_\alpha)| + 7)\\
      &>& -12.
    \end{eqnarray*}
  
{{\bf(\rmnum{3})} $Y_\alpha$ and $Z_\alpha$ are pure, $U_\alpha$ or $V_\alpha$ is contaminated.} 
     
The results of this case are better than the case in which $Z_\alpha$ is contaminated, due to the same test number but more identified items. So $h - \sum_{T\in C_4\cup C'}|I(T)| > -16$.
  
Hence in {Case~2b3}, we have
  \begin{eqnarray*}
    \sum_{T\in C_4\cup C'} |I(T)| &\le& 1.439\cdot d_{C_4\cup C'}\big(\log \frac{\sum_{T\in C_4\cup C'} n(T)}{d_{C_4\cup C'}} + 1.1242\big)+ 16.
  \end{eqnarray*}
  
By {Case~2b1}, {Case~2b2}, and {Case~2b3}, we have the following result in {Case~2b},
  $$
  \sum_{T\in C_4\cup C'} |I(T)| \le 1.431\cdot d_{C_4\cup C'} \big(\log \frac{\sum_{T\in C_4\cup C'} n(T)}{d_{C_4\cup C'}} + 1.1242\big) + 16.
  $$
  
Accordingly, by combining {Case~2a} and {Case~2b}, we have 
  \begin{eqnarray*}
    \sum_{T\in C_2\cup C_3 \cup C_4} |I(T)| &=& \sum_{T\in C_2\cup C_3\setminus C'}|I(T)| + \sum_{T\in C_4\cup C'}|I(T)|\\
    &\le& 1.431\cdot(d_{Z^u} - d_{C_4\cup C'})\cdot\big(\log \frac{\sum_{T\in C_2\cup C_3\setminus C'}n(T)}{d_{Z^u} - d_{C_4\cup C'}}+1.1242\big)\\
    & & +1.431\cdot d_{C_4\cup C'} \big(\log \frac{\sum_{T\in C_4 \cup C'}n(T)}{d_{C_4\cup C'}} + 1.1242\big) + 16\\
    &\le& 1.431\cdot d_{Z^u}\big(\log \frac{\sum_{T\in C_2\cup C_3\cup C_4}n(T)}{d_{Z^u}} + 1.1242\big) + 16.
  \end{eqnarray*}
  
\
  
\noindent{{\bf $\bullet$~Case 3.}}~$|C_4|\ge 1$ and $^1T^1 \in C_4$.
  
Since $^1T^1$ is in $C_4$, one defective item can be detected from some tests in $C_4$. For $C_4$, $\sum_{T\in C_4}|I(T)|$ is equal to $i_\gamma + 1$ plus $|I(^1T^1)|$, and the number of identified items is at least $3\cdot 2^{i_\gamma - 1}$ together with one defective item from $^1T^1$. Since by Claim~\ref{ob:rank}, $r(^1T^1)+ 2 \le r(T^{i_m}_m)$, then, $i_m \ge 3$. Consider two possible cases for $i_\alpha$, i.e., Case 3a and Case 3b.
  
{\noindent{\bf -~Case 3a.}}~$i_\alpha \le 2$.
  
Set $C' = $$\{T^{0}_\alpha, P_m\}$ or $\{P_\alpha, P_m\}$. So three defective items from $D(T^0_\alpha)$ or $D(P_\alpha)$, $P_m$, and $D(^1T^1)$ are identified. And the fact that $I(T^0_\alpha)$ or $I(P_\alpha)$ is no more than $7$ implies that the result of this case is consistent with $i_\beta\le 2$ in {Case~2b2}. Hence, we have
  $$
  \sum_{T\in C_4\cup C'} |I(T)| \le 1.431\cdot\big(\log \frac{\sum_{T\in C_4\cup C'} n(T)}{d_{C_4\cup C'}} + 1.1242\big) + 9.
  $$
  
{\noindent{\bf -~Case 3b.}}~$i_\alpha \ge 3$. 
  
Set $C' = \{P_\alpha, P_m\}$. So three defective items from $P_\alpha, P_m$, and $D(^1T^1)$ are identified. Since the result of this case is consistent with the result of {Case~2b2} where $i_\beta \ge 3$, we have
  $$
  \sum_{T\in C_4\cup C'} |I(T)| \le 1.431\cdot\big(\log \frac{\sum_{T\in C_4\cup C'} n(T)}{d_{C_4\cup C'}} + 1.1242\big) + 10.
  $$
  
Accordingly, by combining {Case~3a} and {Case~3b}, we have
  \begin{eqnarray*}
    \sum_{T\in C_2\cup C_3 \cup C_4} |I(T)| &=& \sum_{T\in C_2\cup C_3\setminus C'}|I(T)| + \sum_{T\in C_4\cup C'}|I(T)|\\
    &\le& 1.431\cdot(d_{Z^u} - d_{C_4\cup C'})\cdot\big(\log \frac{\sum_{T\in C_2\cup C_3\setminus C'}n(T)}{d_{Z^u} - d_{C_4\cup C'}}+1.1242\big)\\
    & & +1.431\cdot d_{C_4\cup C'} \big(\log \frac{\sum_{T\in C_4 \cup C'}n(T)}{d_{C_4\cup C'}} + 1.1242\big) + 10\\
    &\le& 1.431\cdot d_{Z^u}\big(\log \frac{\sum_{T\in C_2\cup C_3\cup C_4}n(T)}{d_{Z^u}} + 1.1242\big) + 10.
  \end{eqnarray*}
  
\
  
Finally, by combining the above three cases {Case~1}, {Case~2}, and {Case~3}, it follows that
  $$
  \sum_{T\in C_2\cup C_3\cup C_4} |I(T)| \le 1.431 d_{Z^u} \big(\log \frac{\sum_{T\in C_2\cup C_3\cup C_4} n(T)}{d_{Z^u}}+1.1242\big) + 16.
  $$
To sum up, Lemma~\ref{C2C3C4} is established. \hfill $\Box$

\end{appendix}

\end{document}